\RequirePackage{fix-cm}
\documentclass[smallextended]{svjour3}       % onecolumn (second format)
\smartqed  % flush right qed marks, e.g. at end of proof
\usepackage{tabularx}
\usepackage{amsmath}
\usepackage{amssymb}
\usepackage{dsfont}
\usepackage{rotating}
\usepackage{multirow}
\usepackage{epstopdf}
\usepackage{newtxmath}
\usepackage{mathtools}

\usepackage{multirow} % tablas con celdas de varias lineas
\usepackage{longtable,tabu} % tablas
\usepackage{hhline}
\usepackage[linesnumbered,ruled,vlined]{algorithm2e}

\usepackage{hyperref} 
\hypersetup{
	colorlinks=true,                          
	linkcolor=black, % Couleur des liens internes
	citecolor=black, % Couleur des num�ros de la biblio dans le corps
	urlcolor=blue  } % Couleur des url
\renewcommand{\vv}{\mathbf{v}}

\newcommand{\q}{\mathbf{q}}

\newcommand{\dt}{{\Delta t}}
\newcommand{\dx}{{\Delta x}}
\newcommand{\dy}{{\Delta y}}
\newcommand{\dz}{{\Delta z}}

\newcommand{\F}{\mathbf{F}}
\newcommand{\G}{\mathbf{G}}
\newcommand{\HH}{\mathbf{H}}
\renewcommand{\S}{\mathbf{S}}

\newcommand{\E}{\mathbf{E}}

\newcommand{\B}{\mathbf{B}}

\newcommand{\g}{\mathbf{g}}

\renewcommand{\vv}{\mathbf{v}}
\newcommand{\fpi}{4\pi}
\newcommand{\epi}{8\pi}
\newcommand{\pdd}[2]{\frac{\partial #1}{\partial #2}}

\newcommand{\Dq}{\Delta \q}
\newcommand{\Dr}{\Delta \rho}
\newcommand{\Dru}{\Delta \rho u}
\newcommand{\Drv}{\Delta \rho v}
\newcommand{\Drw}{\Delta \rho w}
\newcommand{\DrE}{\Delta \rho E}
\newcommand{\Dre}{\Delta \rho e}
\newcommand{\Dp}{\Delta p}
\newcommand{\Dh}{\Delta h}
\newcommand{\DB}{\Delta \B}
\newcommand{\DBx}{\Delta B_x}
\newcommand{\DBy}{\Delta B_y}
\newcommand{\DBz}{\Delta B_z}
\newcommand{\Dm}{\Delta m}

\newcommand{\qt}{\tilde{\q}}
\newcommand{\rt}{\tilde{\rho}}
\newcommand{\ut}{\tilde{u}}
\newcommand{\vt}{\tilde{v}}
\newcommand{\wt}{\tilde{w}}
\newcommand{\pt}{\tilde{p}}
\newcommand{\Bt}{\tilde{B}}
\newcommand{\mt}{\tilde{m}}
\newcommand{\kt}{\tilde{k}}
\newcommand{\qb}{\bar{\q}}
\newcommand{\rb}{\bar{\rho}}
\newcommand{\ub}{\bar{u}}
\newcommand{\vb}{\bar{v}}
\newcommand{\wb}{\bar{w}}

\newcommand{\Bb}{\bar{B}}
\newcommand{\mb}{\bar{m}}
\newcommand{\kb}{\bar{k}}

\newcommand{\RIcolor}[1]{{\color{red}}}
\newcommand{\RIIcolor}[1]{{\color{blue}}}

\newfont{\numerikEleven}{ecrm1000}
\newfont{\numerikTen}{cmss10}
\newfont{\numerikNine}{cmss9}
\newfont{\numerikEight}{cmss8}

\journalname{Journal of Scientific Computing}
\begin{document}
	
\title{A well-balanced semi-implicit IMEX finite volume scheme for ideal Magnetohydrodynamics at all Mach numbers}
	
\titlerunning{A well-balanced semi-implicit finite volume scheme for MHD}        % if too long for running head
	
\author{Claudius Birke \and
		    Walter Boscheri \and
		    Christian Klingenberg
}
	
	
\institute{C. Birke \at
		Department of Mathematics, University of W\"urzburg, Germany \\
		\email{claudius.birke@mathematik.uni-wuerzburg.de}
		\and 
		W. Boscheri \at
		Department of Mathematics and Computer Science, University of
		Ferrara, Italy \\
		\email{walter.boscheri@unife.it}
		\and
		C. Klingenberg \at
		Department of Mathematics, University of W\"urzburg, Germany \\
		\email{klingen@mathematik.uni-wuerzburg.de}
}
	
\date{Received: date / Accepted: date}

\maketitle
	
\begin{abstract}
We propose a second-order accurate semi-implicit and well-balanced finite volume scheme for the equations of ideal magnetohydrodynamics (MHD) including gravitational source terms. 
The scheme treats all terms associated with the acoustic pressure implicitly while keeping the remaining terms part of the explicit sub-system. This semi-implicit approach makes the method particularly well suited for problems in the low Mach regime. We combine the semi-implicit scheme with the deviation well-balancing technique and prove that it maintains equilibrium solutions for the magnetohydrostatic case up to rounding errors. In order to preserve the divergence-free property of the magnetic field enforced by the solenoidal constraint, we incorporate a constrained transport method in the semi-implicit framework. Second order of accuracy is achieved by means of a standard spatial reconstruction technique with total variation diminishing (TVD) property, and by an asymptotic preserving (AP) time stepping algorithm built upon the implicit-explicit (IMEX) Runge-Kutta time integrators. Numerical tests in the low Mach regime and near magnetohydrostatic equilibria support the low Mach and well-balanced properties of the numerical method.

\keywords{semi-implicit \and 
			      well-balanced \and
			      asymptotic-preserving \and
			      compressible low Mach number flows \and
			      finite volume schemes \and
			      magnetohydrodynamics
}
\subclass{MSC 65 \and MSC 68}
\end{abstract}
	
%\tableofcontents
	
%--------- SECTION --------------------------------------------------------
\section{Introduction} \label{sec.intro}

%\cite{BosPar2021}

%MHD
The equations of magnetohydrodynamics (MHD) describe electrically conducting fluids in the presence of a magnetic field. They can be used to model a wide range of physical phenomena, including astrophysical systems such as stars and galaxies, plasma physics, and fusion reactors. 
Their simplest form is given by the ideal MHD equations, in which the influence of fluid viscosity is neglected. This model consists of a system of nonlinear hyperbolic partial differential equations (PDEs) that involve the conservation of mass, momentum, and total energy, along with Faraday's law for the magnetic field. \\
%Low Mach
In certain applications, such as flows in the deep convective layers of stars, the acoustic Mach number, which is the ratio of the fluid speed to the speed of sound, can become very small. Small Mach numbers pose a challenge for explicit finite volume methods, which are typically used to simulate compressible flows, since in this case the flow approaches the incompressible regime and the acoustic waves propagate much faster than the fluid motion. This gives rise to the following problems for standard finite volume methods: the numerical dissipation increases significantly and the time step condition becomes more restrictive. Basically, there are two strategies in the literature to tackle these challenges.
\begin{enumerate}
\item One approach is to rescale the dissipation term in the numerical flux to make it independent of the Mach number and combine this with a fully implicit time integration to circumvent the restrictive time step condition \cite{Viallet2011,BK2022}. The drawback of this strategy is that solving large nonlinear systems of equations is required due to the fully implicit method, which can be computationally expensive.%, and the convergence is difficult to control numerically. 
\item Another approach is to split the PDE system into a convective part that is solved explicitly and an acoustic part that is solved implicitly \cite{SIMHD_Dumbser2019,BosPar2021,Chen2022}. This makes the CFL condition only coupled to the convective sub-system and therefore independent of the Mach number. At the same time, the implicit part becomes smaller and easier to invert, making it less computationally expensive to solve. Furthermore, no numerical dissipation is embedded in the implicit solver.
\end{enumerate}
For our scheme, we pursue the latter approach and treat all terms associated with the acoustic pressure implicitly. All remaining terms are part of the explicit sub-system. For this purpose, a time integration inspired by the class of implicit-explicit (IMEX) schemes is used \cite{BosFil2016,BosPar2021}. Concerning the development of semi-implicit numerical schemes applied to compressible magnetized plasma flows governed by the ideal MHD equations, only very little work has
been done so far. Indeed, semi-implicit schemes for MHD flows have been firstly designed in \cite{SIMHD_Luciani91,SIMHD_Luciani99,SIMHD_Harned86}, however they are concerned with the incompressible case. In \cite{SIMHD_Smolarkiewicz13} the anelastic case of the MHD equations is discussed, while in \cite{SIMHD_Dumbser2019} a semi-implicit finite difference scheme is propsed, inspired by the flux splitting approach originally forwarded for compressible gas dynamics \cite{Dumbser_Casulli16}. An implicit treatment of the magnetic field has been recently devised in \cite{3splitMHD}, and in \cite{Chen2022} a semi-implicit IMEX numerical method for the ideal MHD has been proposed with a finite difference spatial discretization. Our work is strictly related to what presented in \cite{Chen2022}, however there are some important differences: i) we use a finite volume discretization for the convective terms; ii) no numerical dissipation is added in the implicit part even in the case of shock waves; iii) we also include gravitational source terms. \\
%Gravity
Regarding the discretisation in the case of low Mach numbers, this has no effect and the source terms can easily be assigned to the explicit part.
More challenging is its impact on magnetohydrostatic solutions. These are stationary solutions, where the velocity of the fluid is zero and the magnetic field and pressure are in balance with the gravitational force.
If one wants to resolve small perturbations of this equilibrium on a coarse grid, the method should be able to maintain the equilibrium exactly. Usually, this is prevented by different discretizations of the flux and the source term. Numerical schemes that are still able to maintain this balance are called \textit{well-balanced}. There are different strategies in the literature to construct well-balanced schemes, e.g. hydrostatic reconstructions \cite{Kappeli2014,Kappeli2016,Varma2019}, path conservative methods \cite{Pares2006,Castro2008} or special approximate Riemann solvers based on relaxation systems \cite{Desveaux2016,Thomann2020,Birke2023}.
In this work, we combine the semi-implicit scheme with the method of deviation well-balancing \cite{Berberich2021}. Instead of evolving the actual solution, the deviation from an a priori known equilibrium solution is evolved. Through special discretization, the numerical fluxes and source terms cancel out in the case of equilibrium, so that the discrete solution is preserved up to rounding errors. A similar approach has been recently forwarded in \cite{WB_Dumbser2021} in the context of general relativity.\\
%Solenoidal Constraint
In simulations of the MHD equations, it is important to consider the solenoidal constraint, which states that the magnetic field is divergence-free. Numerically, we achieve this structure preserving property by incorporating a constrained transport (CT) method \cite{Gardiner2005} into the semi-implicit scheme. The CT method corrects the magnetic field directly after the explicit step, as the induction equations are entirely assigned to the explicit sub-system of the splitting. By correcting the magnetic field, a specific discrete definition of divergence is kept within machine precision, thereby increasing the stability of the overall scheme. \\

%structure
The paper is structured as follows. In Sect. \ref{sec.pde}, we provide an overview of the MHD equations and their associated properties and propose a splitting into convective and pressure sub-systems. On the basis of this splitting, we develop in Sect. \ref{sec.numscheme} a second order accurate semi-implicit and well-balanced numerical scheme. The numerical results of this scheme for different test problems are presented in Sect. \ref{sec.numtest}, and finally, Sect. \ref{sec.concl} concludes the article providing an outlook to future developments.

%--------- SECTION --------------------------------------------------------
\section{Governing equations} \label{sec.pde}
Let us consider a three-dimensional computational domain $\Omega(x,y,z) \in \mathds{R}^3$ and let the time coordinate be denoted with $t \in \mathds{R}_0^+$. The ideal MHD equations with a gravitational source term constitute a system of balance laws of the form
\begin{equation}
\label{eqn.mhd}
\frac{\partial}{\partial t} 
\left( 
\begin{array}{c}
\rho \\ \rho \vv \\ \rho E \\ \B
\end{array} 
\right)
+
\nabla \cdot
\left( \begin{array}{c}
	\rho \vv \\ 
    \rho \vv \otimes \vv + (p+\frac{1}{2}\vert \B \vert^2) \mathbf{I} - \B \otimes \B \\
    ( \rho E + p + \frac{1}{2}\vert \B \vert^2 ) \vv - \B \left(\B \cdot \vv \right) \\
    \vv \otimes \B - \B \otimes \vv
\end{array} \right)
=
\left( 
\begin{array}{c}
0 \\ \rho \g \\ \rho \vv \cdot \g \\ \mathbf{0}
\end{array} 
\right).
\end{equation}
In this system, $\rho$ denotes the density, $\vv=(u,v,w)$ the velocity field, $p$ the pressure and $\B=(B_x, B_y,B_z)$ the magnetic field. The total energy $\rho E$ is defined as the sum of internal ($\rho e$), kinetic ($\rho k$) and magnetic energy ($m$), i.e.
\begin{equation}
\label{eqn.totEn}
\rho E = \rho e + \rho k + m, \qquad \rho e = \frac{p}{\gamma-1}, \quad \rho k = \frac{1}{2}\rho \vv^2, \quad m = \frac{1}{\epi} \B^2.
\end{equation} 
The internal energy $\rho e$ is computed from the ideal gas equation of state (EOS) with the ratio of specific heats given by $\gamma$. The vector $\g = (g_x,g_y,g_z)$ in the source term on the right hand side of \eqref{eqn.mhd} denotes the time-independent gravitational acceleration.
Using these defintions we can write the MHD equations in the equivalent form
\begin{equation}
	\label{eqn.pde}
	\frac{\partial \q}{\partial t} + \frac{\partial \F(\q)}{\partial x} + \frac{\partial \G(\q)}{\partial y} + \frac{\partial \HH(\q)}{\partial z} = \S(\q),
\end{equation}
with the state vector $\q$, the flux in $x$-direction $\F(\q)$ and the source $\S(\q)$ that explicitly write
\begin{equation}
\label{eqn.mhd_vectors}
\q = \left( \begin{array}{c}
		\rho \\ \rho u \\ \rho v \\ \rho w \\ \rho E \\ B_x \\ B_y \\ B_z
	\end{array} \right), \qquad
\F(\q) = \left( \begin{array}{c}
	\rho u \\ \rho u^2 + p + \frac{1}{\epi} \B^2 - \frac{1}{\fpi} B_x^2 \\ \rho uv - \frac{1}{\fpi} B_x B_y \\ \rho uw - \frac{1}{\fpi} B_x B_z \\ u(\rho E + p + \frac{1}{\epi} \B^2) - \frac{1}{\fpi} B_x (\vv \cdot \B) \\ 0 \\ u B_y - v B_x \\ u B_z - w B_x
\end{array} \right), \qquad
\S(\q)  = \left( \begin{array}{c}
		0 \\ \rho g_x \\ \rho g_y \\ \rho g_z \\ \rho \vv \cdot \g \\ 0 \\ 0 \\ 0
	\end{array} \right).
\end{equation}
The fluxes $\G(\q)$ and $\HH(\q)$ can be expressed in similar forms.

\subsection{Eigenstructure of the MHD system}
In order to analyze the eigenstructure of the MHD system, we consider its one-dimensional homogeneous version. The resulting system is hyperbolic since the eigenvalues $\lambda^{MHD}_{i=\{1,\ldots,8\}}$ of the associated Jacobian matrix $\mathbf{A}=\partial \F/\partial \q$ with $B_x=const$ are
\begin{equation}
	\lambda^{MHD}_{1,8} = u \pm c_f, \quad \lambda^{MHD}_{2,7} = u \pm c_a, \quad \lambda^{MHD}_{3,6} = u \pm c_s, \quad  \lambda^{MHD}_4=u, \quad \lambda^{MHD}_5 = 0.
\end{equation}
The Alfv\'en wave speed ($c_a$) and the slow ($c_s$) and fast ($c_f$) magnetosonic wave speeds herein are defined by
\begin{eqnarray}
	\label{eqn.wavespeeds}
	c_a &=& \frac{B_x}{\sqrt{\fpi \rho}}, \nonumber \\
	c_s^2&=&\frac{1}{2} \left(b^2+c^2-\sqrt{(b+c)^2-4 c_a^2 c^2}\right), \\
	c_f^2&=&\frac{1}{2} \left(b^2+c^2+\sqrt{(b+c)^2-4 c_a^2 c^2}\right), \nonumber
\end{eqnarray}
where the adiabatic sound speed $c= \sqrt{\gamma p / \rho}$ is computed from the ideal gas equation of state. Furthermore, we use the abbreviation $b^2=\B^2/(\fpi \rho)$.

\subsection{Low Mach number regime}
The hydrodynamics behavior of the fluid can be studied by considering the Mach number $\mathcal{M}=\vert \vv \vert/c$. In the low Mach
number limit, the sound speed becomes very high compared to the fluid velocity, hence the terms related to the pressure are dominant. Consequently, larger values of the fast and slow magnetosonic wave speeds are retrieved, and fully explicit numerical methods suffer from both an excessive amount of numerical viscosity, which is proportional to the eigenvalues, and a drastic reduction of the admissible time step $\dt$ to ensure stability under a classical CFL condition of the type
\begin{equation}
	\label{eqn.timestep_ex}
	\dt \leq \text{CFL} \min_{\Omega} \frac{\min (\dx,\dy,\dz)}{\max |\lambda^{MHD}|},
\end{equation}
with $(\dx,\dy,\dz)$ denoting the characteristic mesh spacings in the three spatial directions and the $\text{CFL}\leq1$ being the CFL number. Therefore, we propose to discretize the pressure gradient in the momentum equation and the enthalpy term in the energy equation implicitly, while keeping an explicit discretization for the nonlinear convective fluxes and the terms related to the magnetic field. To that aim, let the fluxes in $x$-direction be split into a convective-type flux $\F^c(\q)$ and a pressure-type flux $\F^p(\q)$, that is
\begin{equation}
	\label{eqn.fluxsplit}
	\F^c(\q) = \left( \begin{array}{c}
		\rho u \\ \rho u^2 + m - \frac{1}{\fpi} B_x^2 \\ \rho uv - \frac{1}{\fpi} B_x B_y \\ \rho uw - \frac{1}{\fpi} B_x B_z \\ u(\rho k + m) - \frac{1}{\fpi} B_x (\vv \cdot \B) \\ 0 \\ u B_y - v B_x \\ u B_z - w B_x
	\end{array} \right), \qquad \F^p(\q) = \left( \begin{array}{c}
	0 \\ p \\ 0 \\ 0 \\ h \rho u \\ 0 \\ 0 \\ 0  
\end{array} \right).
\end{equation} 
We obtain two sub-systems with the following eigenvalues.
\begin{itemize}
	\item Convective sub-system:
	\begin{subequations}
		\label{eqn.convsyst}
		\begin{align}
		&\frac{\partial \q}{\partial t} + \frac{\partial \F^c}{\partial x} = \S, \\
		&\lambda^c_{1,8} = u \pm \sqrt{\frac{\B^2}{\fpi}}, \quad \lambda^c_{2,7} = u \pm \frac{B_x}{\sqrt{\fpi}}, \quad \lambda^c_{3,4} = 0.  	
		\end{align}
	\end{subequations}
    \item Pressure sub-system:
    \begin{subequations}
    	\label{eqn.prssyst}
    	\begin{align}
    	&\frac{\partial \q}{\partial t} + \frac{\partial \F^p}{\partial x} = \mathbf{0}, \\
    	&\lambda^p_{1} = \frac{1}{2} \left( u - \sqrt{u^2+4c^2} \right), \quad \lambda^p_{2,3,4,5,6,7}=0, \quad \lambda^p_8=\frac{1}{2} \left( u + \sqrt{u^2+4c^2} \right).  	
    	\end{align}
    \end{subequations}
\end{itemize}
The fluxes $\G$ and $\HH$ can be split in a similar way with analogous sub-systems and eigenvalues. In the above splitting we have added the source term to the convective sub-system, since this has no effect on the eigenvalues and thus it does not pose a numerical problem in the low Mach limit. It is clear that, by taking the pressure sub-system implicitly, the maximum admissible time step of the scheme becomes
\begin{equation}
	\label{eqn.timestep_im}
	\dt \leq \text{CFL} \min_{\Omega} \frac{\max |\lambda^{c}|}{\min (\dx,\dy,\dz)},
\end{equation}
hence making the scheme particularly well suited for low Mach number flows ($\mathcal{M} \ll 1$) where the pressure terms are dominant. On the other hand, for strongly convected flows with shocks, the convective eigenvalues lead the computation of the time step granting stability.

\subsection{Divergence constraint}
\label{sec.divB}
The induction equations, which describe the evolution of the magnetic field in  \eqref{eqn.mhd}, can be rewritten in the equivalent form 
\begin{equation}
\label{eqn.induction}
\pdd{}{t} \B + \nabla \times \E = 0,
\end{equation}
where $\E= - \vv \times \B$ denotes the electric field. Applying the divergence operator to \eqref{eqn.induction} yields
\begin{equation}
\pdd{}{t} \left( \nabla \cdot \B \right) = 0.
\end{equation}
This implies that the divergence of the magnetic field is always equal to zero if it is zero for the initial data. Therefore, the MHD equations in multiple spatial dimensions naturally contain the condition
\begin{equation}
\label{eq:divB_constraint}
\nabla \cdot \B = 0,
\end{equation}
which is a physically meaningful condition, as it expresses that there should be no magnetic monopoles. \\
At the discrete level, on the other hand, the divergence of the numerically computed magnetic field can increase in time. This can lead to severe stability problems \cite{Brackbill1980}. Therefore, particular care must be taken in the discretisation to ensure that the divergence in the numerical solution remains sufficiently small \cite{BalsaraSpicer1999,divcurlB}. In this work, we rely on a constrained transport method that keeps the divergence at the level of machine accuracy \cite{Gardiner2005}. Sect. \ref{sec.ct} explains how this method works and how it is incorporated into the overall method.

\subsection{Magnetohydrostatic solutions}
\label{sec.mhse}
When considering the MHD equations with gravitational source term, magnethydrostatic solutions can play an important role. These solutions fulfil the following two conditions
\begin{subequations}
\label{eqn.mhse}
\begin{align}
\vv &= 0, \\
\nabla \cdot \left(p+\frac{1}{\epi} \vert \B \vert^2 - \frac{1}{\fpi} \B \otimes \B \right) &= \rho \g.
\end{align}
\end{subequations}
As a consequence, for magnetohydrostatic solutions all fluxes except the momentum flux drop out. The momentum flux, however, fulfils an equilibrium with the source term, so that magnetohydrostatic solutions are stationary.\\
Different ways of discretising the hyperbolic flux and the graviatational source term prevent in general the exact conservation of discrete stationary magnetohydrostatic solutions of the form \eqref{eqn.mhse}.
In order to ensure that the numerical method described in this work is nevertheless exact (up to machine precision) on this type of stationary solutions, we resort to the idea of deviation well-balancing \cite{Berberich2021}. 
The key idea of this approach is to evolve the deviation from an \textit{a priori} known equilibrium solution $\tilde{\q}$ instead of the actual solution $\q$. Therefore, we do not directly discretize equation \eqref{eqn.pde}, but the following equation for the deviation $\Delta \q = \q - \tilde{\q}$:
\begin{equation}
\label{eqn.dq_law}
\begin{split}
&\frac{\partial \Delta \q}{\partial t} 
+ \frac{\partial \F(\Delta \q + \tilde{\q})}{\partial x} - \frac{\partial \F(\tilde{\q})}{\partial x}  
+ \frac{\partial \G(\Delta \q + \tilde{\q})}{\partial y} - \frac{\partial \G(\tilde{\q})}{\partial y} 
+ \frac{\partial \HH(\Delta \q + \tilde{\q})}{\partial z} - \frac{\partial \HH(\tilde{\q})}{\partial z}\\
&= \S(\Delta \q + \tilde{\q}) - \S (\tilde{\q}).
\end{split}
\end{equation}
In the equilibrium case ($\q=\tilde{\q}$), the spatial discretizations of \eqref{eqn.dq_law} cancel so that the equation reduces to $\frac{\partial \Delta \q}{\partial t}=\mathbf{0}$ and the solution remains stationary.

%--------- SECTION --------------------------------------------------------
\section{Numerical scheme} 
\label{sec.numscheme}

For reasons of clarity and comprehensibility we restrict ourselves in the following to one spatial dimension before returning to three spatial dimensions in Section \ref{sec.multi-d}. The computational domain $\Omega = [x_L;x_R]$ is discretized using a total number of $N_x$ equidistant cells of volume $\dx=(x_R-x_L)/N_x$. The cell centers are indicated with $x_i$ and the cell interfaces are referred to with $x_{i+1/2}$. The time coordinate is bounded in the interval $t\in[0;t_f]$, and the final time $t_f$ is reached performing a sequence of time steps $\dt=t^{n+1}-t^n$ that are computed according to the one-dimensional version of the CFL stability condition \eqref{eqn.timestep_im}. 

\subsection{First order semi-discrete scheme in time}
\label{sec.semi-discrete_scheme}
The method constructed in this paper should be able to preserve magnetohydrostatic solutions of the form \eqref{eqn.mhse}. Therefore, we assume to know an equilibrium solution $\tilde{\q}$ \textit{a priori} and discretize equation \eqref{eqn.dq_law}.
The splitting \eqref{eqn.fluxsplit} is therefore applied to equation \eqref{eqn.dq_law}. To simplify the notation, we write $\bar{q}^n = \Delta q^n + \tilde{q}$. For the convective sub-system, which is treated explicitly, we then get:
\begin{equation}
\label{eqn.q-ex}
\Dq^* = \Dq^n - \dt \pdd{}{x} \F^c(\qb) + \dt \pdd{}{x}  \F^c(\qt) + \dt  \S(\qb) - \dt \S(\qt).
\end{equation}
Written out for the single components, this results in
\begin{subequations}
	\begin{align}
		\Dr^{*}    =& \Dr^{n} - \dt \pdd{}{x}\left( \rb \ub \right)^n + \dt \pdd{}{x} (\rt \ut), \label{eqn.rho-ex} \\
		\begin{split}
		(\Dru)^{*} =& (\Dru)^{n} - \dt \pdd{}{x} \left( \rb \ub^2 + \mb - \frac{1}{\fpi} \Bb_x^2 \right)^n + \dt \pdd{}{x} \left( \rt \ut^2 + \mt - \frac{1}{\fpi} \Bt_x^2 \right) \\
		&+ \dt \left( \rb g_x \right) - \dt \left( \rt g_x \right), \label{eqn.rhou-ex}
		\end{split} \\
		(\Drv)^{*}=& (\Drv)^{n} - \dt \pdd{}{x} \left( \rb \ub \vb - \frac{1}{\fpi} \Bb_x \Bb_y \right)^n + \dt \pdd{}{x} \left( \rt \ut \vt - \frac{1}{\fpi} \Bt_x \Bt_y \right), \label{eqn.rhov-ex} \\
		(\Drw)^{*}=& (\Drw)^{n} - \dt \pdd{}{x} \left( \rb \ub \wb - \frac{1}{\fpi} \Bb_x \Bb_z \right)^n + \dt \pdd{}{x} \left( \rt \ut \wt - \frac{1}{\fpi} \Bt_x \Bt_z \right), \label{eqn.rhow-ex} \\
		\begin{split}
		(\DrE)^{*} =& (\DrE)^{n} - \dt \pdd{}{x} \left( \ub (\rb \kb + \mb) - \frac{1}{\fpi} \Bb_x (\bar{\vv} \cdot \bar{\B})  \right)^n + \dt \pdd{}{x} \left( \ut (\rt \kt + \mt) - \frac{1}{\fpi} \Bt_x (\tilde{\vv} \cdot \tilde{\B}) \right) \\
		&+ \dt \left(\rb \ub g_x \right) - \dt \left(\rt \ut g_x \right), \label{eqn.rhoE-ex}
		\end{split} \\
		\DBx^{*}     =& 0, \label{eqn.Bx-ex} \\
		\DBy^{*}     =& \DBy^{n} - \dt \pdd{}{x} \left( \ub \Bb_y - \vb \Bb_x \right)^n + \dt \pdd{}{x} \left( \ut \Bt_y - \vt \Bt_x \right), \label{eqn.By-ex} \\
		\DBz^{*}     =& \DBz^{n} - \dt \pdd{}{x} \left( \ub \Bb_z - \wb \Bb_x \right)^n + \dt \pdd{}{x} \left( \ut \Bt_z - \wt \Bt_x \right). \label{eqn.Bz-ex}
	\end{align}
\end{subequations}
%\begin{subequations}
%	\begin{align}
%		\Dr^{*}    &= \Dr^{n} - \dt \pdd{}{x}\left( \rb \ub \right)^n, \label{eqn.rho-ex} \\
%		(\Dru)^{*}&= (\Dru)^{n} - \dt \pdd{}{x} \left( \rb \ub^2 + \mb - \frac{1}{\fpi} \Bb_x^2 \right)^n + \dt \pdd{}{x} \left(\mt - \frac{1}{\fpi} \Bt_x^2 \right) + \dt \left( \rb g_x \right) - \dt \left( \rt g_x \right), \label{eqn.rhou-ex} \\
%		(\Drv)^{*}&= (\Drv)^{n} - \dt \pdd{}{x} \left( \rb \ub \vb - \frac{1}{\fpi} \Bb_x \Bb_y \right)^n + \dt \pdd{}{x} \left( - \frac{1}{\fpi} \Bt_x \Bt_y \right), \label{eqn.rhov-ex} \\
%		(\Drw)^{*}&= (\Drw)^{n} - \dt \pdd{}{x} \left( \rb \ub \wb - \frac{1}{\fpi} \Bb_x \Bb_z \right)^n + \dt \pdd{}{x} \left( - \frac{1}{\fpi} \Bt_x \Bt_z \right), \label{eqn.rhow-ex} \\
%		(\DrE)^{*}&= (\DrE)^{n} - \dt \pdd{}{x} \left( \ub (\rb \kb + \mb) - \frac{1}{\fpi} \Bb_x (\bar{\vv} \cdot \bar{\B}) \right)^n  + \dt \left(\rb \ub g_x \right), \label{eqn.rhoE-ex} \\
%		\DBx^{*}     &=0, \label{eqn.Bx-ex} \\
%		\DBy^{*}     &=\DBy^{n} - \dt \pdd{}{x} \left( \ub \Bb_y - \vb \Bb_x \right)^n, \label{eqn.By-ex} \\
%		\DBz^{*}     &=\DBz^{n} - \dt \pdd{}{x} \left( \ub \Bb_z - \wb \Bb_x \right)^n. \label{eqn.Bz-ex}
%	\end{align}
%\end{subequations}
In the above equations, we have included terms containing $\ut$, $\vt$ or $\wt$. In the following, however, we will omit these terms, since the equilibrium velocities are always zero.\\
The above definitions are then employed to obtain a first order semi-implicit time discretization \cite{BosFil2016,BosPar2021,BT2023} of the MHD equations \eqref{eqn.pde}, which writes
\begin{subequations}
	\label{eqn.sd-first}
	\begin{align}
		\Dr^{n+1}    &= \Dr^{*}, \label{eqn.sd-rho} \\
		(\Dru)^{n+1}&=(\Dru)^{*} - \dt \pdd{}{x}\left(p^{n+1} - \pt \right), \label{eqn.sd-rhou} \\
		(\Drv)^{n+1}&=(\Drv)^{*}, \label{eqn.sd-rhov} \\
		(\Drw)^{n+1}&=(\Drw)^{*}, \label{eqn.sd-rhow} \\		 
%		(\Dre)^{n+1}+\frac{(\Dru)^{n+1}\,(\Dru)^{n}}{2 \Dr^{n+1}}+\Delta m^{n+1}&=(\DrE)^{*} - \dt \pdd{}{x} \left( \left(\Dh^n+\tilde{h} \right) \left( (\Dru)^{n+1} + \rt \ut \right) \right), \label{eqn.sd-rhoE} \\
		(\DrE)^{n+1}&=(\DrE)^{*} - \dt \pdd{}{x} \left( h^n (\rho u)^{n+1} \right), \label{eqn.sd-rhoE} \\
		\DBx^{n+1}&=\DBx^{*}, \\
		\DBy^{n+1}&=\DBy^{*}, \\
		\DBz^{n+1}&=\DBz^{*}.
	\end{align}
\end{subequations}
Here, we use the known enthalpy $h^n$ at time level $n$ in the energy equation to avoid nonlinear terms in the implicit part, differently form the schemes proposed in \cite{SIMHD_Dumbser2019,3splitMHD}. Additionally, we have made use of the fact that in the implicit fluxes we only work with point values at the cell centers, i.e. we use an implicit finite difference discretization, thus we can write
\begin{equation}
\left( \Delta \q + \tilde{\q} \right)_i = \q_i - \tilde{\q}_i + \tilde{\q}_i = \q_i.
\end{equation}
Therefore, the implicit fluxes are simplified by
\begin{align*}
&\pdd{}{x}\left(\Dp^{n+1} + \pt \right) - \dt \pdd{}{x} \pt = \pdd{}{x}\left(p^{n+1} - \pt \right), \\
&\pdd{}{x} \left( \left(\Dh^n+\tilde{h} \right) \left( (\Dru)^{n+1} + \rt \ut \right) \right) = \pdd{}{x} \left( h^n (\rho u)^{n+1} \right).
\end{align*}
According to the definition given in the energy equation \eqref{eqn.totEn}, we split the deviation of the total energy at the new time level as follows
\begin{equation}
(\DrE)^{n+1} = (\Dre)^{n+1} + (\Delta m)^{n+1} + (\Delta \rho k)^{n+1}. 
\end{equation}
The deviation of the kinetic energy therein is set to
\begin{equation}
(\Delta \rho k)^{n+1} = (\rho k)^{n+1} - (\rt \tilde{k}) = \frac{1}{2} \frac{(\rho u)^{n}}{\rho^{n+1}} (\rho u)^{n+1} - \frac{1}{2} \rt \ut^2,
\end{equation}
where a semi-implicit strategy is adopted for the term $(\rho k)^{n+1}$. Indeed, the split of the momentum contribution into an explicit and an implicit part is again done in order to avoid nonlinear implicit terms. The density $\rho^{n+1}$ and the deviation in the magnetic energy $\Dm^{n+1}$ are known because both continuity and induction equations are fully explicit. With the help of the ideal gas law, the deviation in the internal energy can be expressed in terms of the pressure
\begin{equation}
(\Dre)^{n+1} = \frac{p^{n+1}-\pt}{\gamma-1}.
\end{equation}
In order to derive a preliminary discretization of the total energy equation we transform \eqref{eqn.sd-rhou} into
\begin{equation}
(\rho u)^{n+1} = (\Dru)^{*} + \rt \ut - \dt \pdd{}{x}\left(p^{n+1} - \pt \right),
\end{equation}
and insert the result into the energy equation \eqref{eqn.sd-rhoE}.
This leads to an elliptic equation for the pressure
\begin{equation}
	\label{eqn.p}
	\frac{p^{n+1}}{\gamma-1} - \dt \frac{(\rho u)^{n}}{2 \rho^{n+1}} \pdd{}{x} \left(p^{n+1}-\pt \right) - \dt^2 \pdd{}{x} \left( h^n \pdd{}{x} \left( p^{n+1} - \pt \right) \right) = b^n,
\end{equation} 
with the known right-hand-side given by
\begin{equation}
	\label{eqn.rhs}
	b^n = \frac{\pt}{\gamma-1} + \frac{1}{2} \rt \ut^2 + (\DrE)^{*} - \frac{(\rho u)^{n}}{2 \rho^{n+1}} \left( (\Dru)^{*} + \rt \ut \right) - \Dm^{n+1} - \dt \pdd{}{x} \left( h^n \left( (\Dru)^{*} + \rt \ut \right) \right).
\end{equation}
%\rev{Comment: We could again drop all parts containing $\ut$.}
The pressure equation \eqref{eqn.p} constitutes a linear system for the scalar unknown $p^{n+1}$ that is solved using the iterative GMRES solver \cite{GMRES} up to a prescribed tolerance (we typically set $\text{tol}=10^{-12}$). Differently from \cite{SIMHD_Dumbser2019,3splitMHD}, this approach does not need any fixed point method thanks to the semi-implicit splitting of the enthalpy flux and the kinetic energy in the energy equation. Once the new pressure is known, the deviation in the momentum $(\Dru)^{n+1}$ is updated with \eqref{eqn.sd-rhou}, and then the deviation in the total energy is updated using \eqref{eqn.sd-rhoE}. Notice that the scheme is written in flux form, therefore it is locally and globally conservative.

\subsection{Discrete spatial operators}
\label{sec.spatial_operators}

The convective sub-system \eqref{eqn.q-ex} is discretized by a Godunov-type finite volume method which writes
\begin{equation}
\label{eqn.spatial_operators-ex}
\Dq_{i}^* = \Dq_{i}^n -\frac{\dt}{\dx} (\hat{\F}_{i+1/2}^{c} - \hat{\F}_{i-1/2}^{c}) + \frac{\dt}{\dx} (\F^c(\qt_{i+1/2}) - \F^c(\qt_{i-1/2})) + \dt \hat{\S}_i -  \dt \S(\qt_i).
\end{equation}
In this notation, the hat symbol indicates a numerical flux or source term. For the numerical flux, we decide to use a simple and robust Rusanov-type flux of the form
\begin{equation}
\label{eqn.rusanov}
\hat{\F}_{i+1/2}^{c} = \hat{\F}(\qb_{i+1/2}^{L},\qb_{i+1/2}^{R}) = \frac{1}{2} \left(\F^c(\qb_{i+1/2}^{L})+\F^c(\qb_{i+1/2}^{R}) \right) - \frac{1}{2} s_{max} \left(\qb_{i+1/2}^{R} - \qb_{i+1/2}^{L} \right),
\end{equation}
where the numerical dissipation $s_{max}=\max \left( \vert (\lambda^c)_{i+1/2}^L \vert,\vert (\lambda^c)_{i+1/2}^R \vert \right)$ only considers the convective eigenvalues and is therefore independent of the Mach number. The numerical flux is computed in the states
\begin{equation}
\label{eqn.interface_states}
\qb_{i+1/2}^{L,R} = \Dq_{i+1/2}^{L,R} + \qt_{i+1/2},
\end{equation}
where the superscripts $L,R$ denote the left and right extrapolated data at the interface. It is essential for the well-balancing that only the deviation $\Dq$ is reconstructed, while the equilibrium solution is evaluated at the cell interface.\\
The third term in the right hand side of \eqref{eqn.spatial_operators-ex} simply computes the physical fluxes based on the equilibrium solution at the cell interface, that is $\F^c(\qb_{i+1/2}^{L,R})$. For the source term, we substitute the volume-averaged quantity with its cell centered value, which is accurate up to second order:
\begin{equation}
\hat{\S}_i = \S(\bar{\q}_i).
\end{equation}
The implicit terms appearing in the pressure sub-system \eqref{eqn.prssyst} are approximated by means of finite difference operators with no numerical dissipation \cite{BosPar2021}:
\begin{subequations}
	\label{eqn.fd}
	\begin{align}
	 \left. \pdd{p}{x} \right|_i^{n+1} &= \frac{p_{i+1}^{n+1}-p_{i-1}^{n+1}}{2 \, \dx} + \mathcal{O}(\dx^2), \label{eqn.fddx} \\
	\pdd{}{x} \left. \left( h \pdd{p}{x} \right) \right|_i^{n,n+1} &= \frac{1}{\dx^2} \left[ h_{i-1}^n \, h_i^n \, h_{i+1}^n \right] \left[ \begin{array}{ccc}
	 	3/4 & -1 & 1/4 \\ 0 & 0 & 0 \\ 1/4 & -1 & 3/4
	 \end{array} \right] \left[ \begin{array}{c}
	 	p_{i-1}^{n+1} \\ p_i^{n+1} \\ p_{i+1}^{n+1}
	 \end{array} \right] + \mathcal{O}(\dx^2). \label{eqn.fdhdx}
	\end{align}
\end{subequations}

\subsection{Second order extension in time and space}
\label{sec.second_order}

\subsubsection{Second order in time}
The scheme achieves second order of accuracy in time by using the implicit-explicit Runge Kutta method LSDIRK(2,2,2) \cite{PR_IMEX}. 
The triplet (2,2,2) refers to the number of stages of the implicit part, the number of stages of the explicit part, and the order of the IMEX scheme, respectively. Let us denote the spatial operator in the MHD equations by $\mathcal{H}(\q_E(t),\q_I(t))$, where the first argument is discretized explicitly, while the second argument is discretized implicitly. According to the splitting \eqref{eqn.fluxsplit}, $\mathcal{H}$ is defined by
\begin{equation}
\mathcal{H}(\q_E(t),\q_I(t)) = \pdd{}{x} \F^c(\q_E) + \pdd{}{x} \F^p(\q_I) + \S(\q_E).
\end{equation}
At the beginning of each time step we initialize $\q_E^n=\q_I^n=\q^n$ and then compute for each stage $i$ its stage fluxes $\mathbf{k}_i$ in the following way:
\begin{subequations}
\begin{align}
\q_E^{i} &= \q_E^n + \dt \sum_{j=1}^{i-1} \hat{a}_{ij} \mathbf{k}_j, &2 \leq i \leq s, \\
\q_I^{i} &= \q_E^n + \dt \sum_{j=1}^{i-1} a_{ij} \mathbf{k}_j, &2 \leq i \leq s, \\
\mathbf{k}_i &= \mathcal{H} \left( \q_E^{i}, \q_I^{i} + \dt a_{ii} \mathbf{k}_i \right), &1 \leq i \leq s. \label{eqn.ki_imex}
\end{align}
\end{subequations}
Notice that the only implicit evaluation lies in \eqref{eqn.ki_imex}, which indeed correspond to the solution of the elliptic equation on the pressure \eqref{eqn.p}.
Using these stage fluxes, the updated solution at the new time level is computed by
\begin{equation}
\q^{n+1} = \q^n + \dt \sum_{i=1}^s b_i \mathbf{k}_i.
\end{equation}
The coefficients are given by a double Butcher tableau of the form
\begin{equation}
\begin{array}
{c|c}
\hat{c} & \hat{A} \\
\hline
 & \hat{b}^\top\\
\end{array}
\qquad
\begin{array}
{c|c}
c & A \\
\hline
 & b^\top\\
\end{array}
\end{equation}
with matrices $(\hat{A},A) \in \mathbb{R}^{s\times s}$ and vectors $(\hat{c},c,\hat{b},b) \in \mathbb{R}^s$. In this notation, the hat symbol indicates the coefficients for the explicit scheme. The coefficients for LSDIRK2(2,2,2) are listed in Table \ref{tab.butcher}.

\begin{table}[h!]
\begin{center}
\begin{tabular}{c|cc}
0           & 0                  & 0 \\
$\beta$ & $\beta$        & 0 \\
\hline
            & 1-$\gamma$ & $\gamma$
\end{tabular}
\quad \quad \quad
\begin{tabular}{c|cc}
$\gamma$ & $\gamma$    & 0 \\
1                & 1-$\gamma$ & $\gamma$ \\
\hline
            & 1-$\gamma$ & $\gamma$
\end{tabular}
\end{center} 
\caption{Butcher tableau for the LSDIRK2(2,2,2) time discretization with  $\gamma = 1-1/\sqrt{2}$, $\beta=1/(2\gamma)$. Left: explicit tableau. Right: implicit tableau.} \label{tab.butcher}
\end{table}

\subsubsection{Second order in space}
The extension to second order in space in the explicit part is based on a piecewise linear reconstruction of the conservative variables. The values at the cell interface, which serve as initial data for the Riemann problems, are computed by evaluating the function
\begin{equation}
\q(x) = \q_i + \mathbf{\sigma}_i (x-x_i)
\end{equation}
at the cell boundaries $x_{i-1/2}$ and $x_{i+1/2}$. The slopes $\mathbf{\sigma}_i$ are computed by applying a minmod limiter to the left and right slope, i.e.
\begin{equation}
\mathbf{\sigma}_i = \text{minmod} \left( \frac{\q_i-\q_{i-1}}{\dx}, \frac{\q_{i+1}-\q_{i}}{\dx} \right).
\end{equation}
We compute a separate slope for each of the conservative variables in $\q$.
At this point, it should be recalled that the reconstruction is only applied to the deviation $\Dq$ in \eqref{eqn.interface_states}, and not to the equilibrium solution $\qt$. 
With this procedure, the spatial discretization of the explicit part becomes second order. In the implicit part, no further changes are needed, since only point values are used and the operators in \eqref{eqn.fd} are already second order.

\subsection{Multi-dimensional extension}
\label{sec.multi-d}
In three spatial dimensions we discretize the domain $\Omega= [x_L,x_R] \times [y_L,y_R] \times [z_L,z_R]$ by a Cartesian grid with $N_x \times N_y \times N_z$ cells, which have the uniform cell size $\dx \times \dy \times \dz$.\\
The discretization of the explicit part is extended to three spatial dimensions by using an unsplit finite volume method according to \cite{ToroBook}, which writes
\begin{equation}
\label{eqn.3d-ex}
\begin{split}
\Dq_{i,j,k}^* =  \Dq_{i,j,k}^n 
&-\frac{\dt}{\dx} \left(\hat{\F}_{i+1/2,j,k}^{c} - \hat{\F}_{i-1/2,j,k}^{c}\right) + \frac{\dt}{\dx} \Bigl( \F^c(\qt_{i+1/2,j,k}) - \F^c(\qt_{i-1/2,j,k}) \Bigr) \\
&-\frac{\dt}{\dy}\left(\hat{\G}_{i,j+1/2,k}^{c} - \hat{\G}_{i,j-1/2,k}^{c}\right) + \frac{\dt}{\dy} \Bigl( \G^c(\qt_{i,j+1/2,k}) - \G^c(\qt_{i,j-1/2,k})\Bigr)  \\
&-\frac{\dt}{\dz} \left(\hat{\HH}_{i,j,k+1/2}^{c} - \hat{\HH}_{i,j,k-1/2}^{c}\right) + \frac{\dt}{\dz} \Bigl(\HH^c(\qt_{i,j,k+1/2}) - \HH^c(\qt_{i,j,k-1/2})\Bigr)  \\
&+ \dt \hat{\S}_{i,j,k} -  \dt \S(\qt_{i,j,k}).
\end{split}
\end{equation}
The numerical fluxes $\hat{\F}^c$, $\hat{\G}^c$ and $\hat{\HH}^c$ have the form of the Rusanov flux and are constructed as in \eqref{eqn.rusanov}. For the source term $\hat{\S}$ we again use the cell centered value. The updated density $\Dr_{i,j,k}^{n+1}$ and magnetic field $\DB_{i,j,k}^{n+1}$ are equal to their explicit update, since in the splitting the complete flux of these components is explicit. The update of the momentum components, on the other hand, contains implicit parts and reads in the fully discrete and three-dimensional form as
\begin{subequations}
\label{eqn.3d-momentum}
\begin{align}
(\Dru)_{i,j,k}^{n+1}&=(\Dru)_{i,j,k}^{*} - \frac{\dt}{2\dx}\left(p_{i+1,j,k}^{n+1} - \pt_{i+1,j,k} - p_{i-1,j,k}^{n+1} + \pt_{i-1,j,k}  \right), \label{eqn.3d-rhou} \\
(\Drv)_{i,j,k}^{n+1}&=(\Drv)_{i,j,k}^{*} - \frac{\dt}{2\dy}\left(p_{i,j+1,k}^{n+1} - \pt_{i,j+1,k} - p_{i,j-1,k}^{n+1} + \pt_{i,j-1,k}  \right), \label{eqn.3d-rhov} \\
(\Drw)_{i,j,k}^{n+1}&=(\Drw)_{i,j,k}^{*} - \frac{\dt}{2\dz}\left(p_{i,j,k+1}^{n+1} - \pt_{i,j,k+1} - p_{i,j,k-1}^{n+1} + \pt_{i,j,k-1}  \right). \label{eqn.3d-rhow}
\end{align}
\end{subequations}
Implicit terms also appear in the update of the total energy:
\begin{equation}
\label{eqn.3d-rhoE}
\begin{split}
(\DrE)_{i,j,k}^{n+1} = (\DrE)_{i,j,k}^* 
&- \frac{\dt}{2\dx} \left( h_{i+1,j,k}^n(\rho u)_{i+1,j,k}^{n+1} - h_{i-1,j,k}^n(\rho u)_{i-1,j,k}^{n+1} \right) \\
&- \frac{\dt}{2\dy} \left( h_{i,j+1,k}^n(\rho v)_{i,j+1,k}^{n+1} - h_{i,j-1,k}^n(\rho v)_{i,j-1,k}^{n+1} \right) \\
&- \frac{\dt}{2\dz} \left( h_{i,j,k+1}^n(\rho w)_{i,j,k+1}^{n+1} - h_{i,j,k-1}^n(\rho w)_{i,j,k-1}^{n+1} \right).
\end{split}
\end{equation}
The pressure $p^{n+1}$, which is needed for the updates \eqref{eqn.3d-momentum} and \eqref{eqn.3d-rhoE}, can be determined by solving the following elliptic equation:
\begin{equation}
\label{eqn.3d-linSyst_p}
\begin{split}
&\frac{p_{i,j,k}^{n+1}}{\gamma-1} 
- \frac{\dt}{2\dx} \frac{(\rho u)_{i,j,k}^{n}}{2 \rho_{i,j,k}^{n+1}} \left(p_{i+1,j,k}^{n+1} - \pt_{i+1,j,k} - p_{i-1,j,k}^{n+1} + \pt_{i-1,j,k}  \right) \\
&- \frac{\dt}{2\dy} \frac{(\rho v)_{i,j,k}^{n}}{2 \rho_{i,j,k}^{n+1}} \left(p_{i,j+1,k}^{n+1} - \pt_{i,j+1,k} - p_{i,j-1,k}^{n+1} + \pt_{i,j-1,k}  \right) \\
&- \frac{\dt}{2\dz} \frac{(\rho w)_{i,j,k}^{n}}{2 \rho_{i,j,k}^{n+1}} \left(p_{i,j,k+1}^{n+1} - \pt_{i,j,k+1} - p_{i,j,k-1}^{n+1} + \pt_{i,j,k-1}  \right) \\
&\begin{split}
- \frac{\dt^2}{\dx^2} \Bigg [ &\left( \frac{3}{4} h_{i-1,j,k}^n + \frac{1}{4} h_{i+1,j,k}^n \right) \left(p_{i-1,j,k}^{n+1}-\pt_{i-1,j,k}\right) 
- \left( h_{i-1,j,k}^n+h_{i+1,j,k}^n \right) \left(p_{i,j,k}^{n+1}-\pt_{i,j,k}\right) \\
&+ \left( \frac{1}{4} h_{i-1,j,k}^n+\frac{3}{4} h_{i+1,j,k}^n \right) \left(p_{i+1,j,k}^{n+1}-\pt_{i+1,j,k} \right) \Bigg ] 
\end{split} \\
&\begin{split}
- \frac{\dt^2}{\dy^2} \Bigg [ &\left( \frac{3}{4} h_{i,j-1,k}^n + \frac{1}{4} h_{i,j+1,k}^n \right) \left(p_{i,j-1,k}^{n+1}-\pt_{i,j-1,k}\right) - \left( h_{i,j-1,k}^n+h_{i,j+1,k}^n \right) \left(p_{i,j,k}^{n+1}-\pt_{i,j,k}\right) \\
&+ \left( \frac{1}{4} h_{i,j-1,k}^n+\frac{3}{4} h_{i,j+1,k}^n \right) \left(p_{i,j+1,k}^{n+1}-\pt_{i,j+1,k}\right) \Bigg] 
\end{split} \\
&\begin{split}
- \frac{\dt^2}{\dz^2} \Bigg [ &\left( \frac{3}{4} h_{i,j,k-1}^n + \frac{1}{4} h_{i,j,k+1}^n \right) \left(p_{i,j,k-1}^{n+1}-\pt_{i,j,k-1}\right) 
- \left( h_{i,j,k-1}^n+h_{i,j,k+1}^n \right) \left(p_{i,j,k}^{n+1}-\pt_{i,j,k}\right) \\
&+ \left( \frac{1}{4} h_{i,j,k-1}^n+\frac{3}{4} h_{i,j,k+1}^n \right) \left(p_{i,j,k+1}^{n+1}-\pt_{i,j,k+1}\right) \Bigg ] 
\end{split} \\
&= b_{i,j,k}^n
\end{split}
\end{equation}

with the right-hand side
%\begin{equation}
%\begin{split}
%b_{i,j,k}^n =& 
%\frac{\pt_{i,j,k}}{\gamma-1} + \frac{1}{2} \rt_{i,j,k} \ut_{i,j,k}^2 + \frac{1}{2} \rt_{i,j,k} \vt_{i,j,k}^2 + \frac{1}{2} \rt_{i,j,k} \wt_{i,j,k}^2 + (\DrE)_{i,j,k}^* - \Dm_{i,j,k}^{n+1} \\
%&-\frac{(\rho u)_{i,j,k}^n}{2\rho_{i,j,k}^{n+1}} \left( (\Dru)_{i,j,k}^* + \rt_{i,j,k} \ut_{i,j,k} \right)
%  -\frac{(\rho v)_{i,j,k}^n}{2\rho_{i,j,k}^{n+1}} \left( (\Drv)_{i,j,k}^* + \rt_{i,j,k} \vt_{i,j,k} \right) \\
% &-\frac{(\rho w)_{i,j,k}^n}{2\rho_{i,j,k}^{n+1}} \left( (\Drw)_{i,j,k}^* + \rt_{i,j,k} \wt_{i,j,k} \right) \\
%&-\frac{\dt}{2\dx} \left( h_{i+1,j,k}^n \left( (\Dru)_{i+1,j,k}^* + \rt_{i+1,j,k} \ut_{i+1,j,k} \right) - h_{i-1,j,k}^n \left( (\Dru)_{i-1,j,k}^* + \rt_{i-1,j,k} \ut_{i-1,j,k} \right) \right) \\
%&-\frac{\dt}{2\dy} \left( h_{i,j+1,k}^n \left( (\Drv)_{i,j+1,k}^* + \rt_{i,j+1,k} \vt_{i,j+1,k} \right) - h_{i,j-1,k}^n \left( (\Drv)_{i,j-1,k}^* + \rt_{i,j-1,k} \vt_{i,j-1,k} \right) \right) \\
%&-\frac{\dt}{2\dz} \left( h_{i,j,k+1}^n \left( (\Drw)_{i,j,k+1}^* + \rt_{i,j,k+1} \ut_{i,j,k+1} \right) - h_{i,j,k-1}^n \left( (\Drw)_{i,j,k-1}^* + \rt_{i,j,k-1} \wt_{i,j,k-1} \right) \right).
%\end{split}
%\end{equation}
\begin{equation}
\begin{split}
b_{i,j,k}^n =& 
\frac{\pt_{i,j,k}}{\gamma-1} + (\DrE)_{i,j,k}^* - \Dm_{i,j,k}^{n+1} \\
&-\frac{(\rho u)_{i,j,k}^n}{2\rho_{i,j,k}^{n+1}} (\Dru)_{i,j,k}^* 
   -\frac{\dt}{2\dx} \left( h_{i+1,j,k}^n (\Dru)_{i+1,j,k}^* - h_{i-1,j,k}^n  (\Dru)_{i-1,j,k}^* \right) \\
& -\frac{(\rho v)_{i,j,k}^n}{2\rho_{i,j,k}^{n+1}} (\Drv)_{i,j,k}^* 
    -\frac{\dt}{2\dy} \left( h_{i,j+1,k}^n (\Drv)_{i,j+1,k}^* - h_{i,j-1,k}^n  (\Drv)_{i,j-1,k}^* \right) \\
&  - \frac{(\rho w)_{i,j,k}^n}{2\rho_{i,j,k}^{n+1}} (\Drw)_{i,j,k}^* 
    -\frac{\dt}{2\dz} \left( h_{i,j,k+1}^n (\Drw)_{i,j,k+1}^* - h_{i,j,k-1}^n (\Drw)_{i,j,k-1}^*  \right).
\end{split}
\end{equation}

\subsection{Constrained transport method}
\label{sec.ct}

When simulating the ideal MHD equations in multiple space dimensions, the divergence constraint \eqref{eq:divB_constraint} must also be taken into account in order to maintain stability and accuracy. Furthermore, the solenoidal property of the magnetic field represents a physical constraint. The scheme described up to this point does not do this in the sense that the divergence can increase significantly during the simulation. Therefore, an additional correction of the magnetic field is required after each time step. In order to retrieve at the discrete level a solenoidal magnetic field, we rely on a second order accurate CT method \cite{Gardiner2005}.
As it is typical in staggered CT methods, this method is based on the integration of the induction equations \eqref{eqn.induction} over the cell boundaries using Stokes theorem. This results in an update formula for the magnetic field at the face centers of the cell based on the electric field at the cell corners. In our case, we update the deviation at the face centers:
\begin{equation}
\label{eqn.Bfc-update}
\begin{split}
\pdd{}{t} \Delta B_{x,i+1/2,j,k} =& 
\frac{1}{\dy} \left( \Delta \E_{z,i+1/2,j+1/2,k} - \Delta \E_{i+1/2,j-1/2,k} \right) \\
&- \frac{1}{\dz} \left( \Delta \E_{y,i+1/2,j,k+1/2} - \Delta \E_{y,i+1/2,j,k-1/2} \right).
\end{split}
\end{equation}
The deviation of the magnetic field at $y$- and $z$-faces is updated by similar formulae. The discretizations of the electric field at the cell corners can be taken from \cite{Gardiner2005} and transferred to the deviations without any problems, since the physical fluxes evaluated in the equilibrium solution are omitted due to the zero velocity. The cell averages of the magnetic field is then computed via the arithmetic mean from the values at the faces, which still maintains second order of accuracy:
\begin{equation}
\Delta B_{x,i,j,k} = \frac{1}{2} \left( \Delta B_{x,i-1/2,j,k} + \Delta B_{x,i+1/2,j,k} \right).
\end{equation}
In the semi-implicit scheme, this whole procedure of constained transport is carried out directly after the explicit update. Thus, for solving the pressure equation \eqref{eqn.3d-linSyst_p} we already use the corrected and divergence-free magnetic field. 
Overall, by incorporating the CT method, the semi-implicit scheme keeps the following discrete definition of the divergence
\begin{equation}
(\nabla \cdot \B)_{i,j,k} = 
\frac{B_{x,i+1/2,j,k}-B_{x,i-1/2,j,k}}{\dx} +
\frac{B_{y,i,j+1/2,k}-B_{y,i,j-1/2,k}}{\dy} +
\frac{B_{z,i,j,k+1/2}-B_{z,i,j,k-1/2}}{\dz},
\end{equation}
within the order of machine accuracy. This is checked and validated at each time step of the computation.

\subsection{Well-balanced property}
\label{sec.wb_proof}

The following theorem states that the presented fully discrete scheme preserves magnetohydrostatic equilibria of the form \eqref{eqn.mhse} up to rounding errors.

\begin{theorem}
\label{th.wb_pressure}
Let us assume that the numerical solution $\q^n$ is equal to the discrete equilibrium solution $\qt$, i.e.
\begin{equation}
\label{eqn.wbp}
\q_{i,j,k}^n = \qt_{i,j,k} \qquad \forall (i,j,k) \in \{1,...,N_x\} \times \{1,...,N_y\} \times \{1,...,N_z\}.
\end{equation}
Then the numerical method described in Sect. \ref{sec.semi-discrete_scheme}-\ref{sec.ct} preserves the solution up to machine precision.
\end{theorem}

\begin{proof}
From the assumption \eqref{eqn.wbp} it follows that $\Dq = 0$ at every grid point. Thus, the input of the Rusanov fluxes in the convective sub-system is reduced to the equilibrium solution at the interface, i.e.
\begin{subequations}
\begin{align}
\qb_{i+1/2,j,k}^{L} = \qb_{i+1/2,j,k}^{R} =  \qt_{i+1/2,j,k}, \\
\qb_{i,j+1/2,k}^{L} = \qb_{i,j+1/2,k}^{R} =  \qt_{i,j+1/2,k}, \\
\qb_{i,j,k+1/2}^{L} = \qb_{i,j,k+1/2}^{R} =  \qt_{i,j,k+1/2}.
\end{align}
\end{subequations}
It follows by the consistency of the Rusanov flux that the numerical fluxes are equal to the physical fluxes and these truncate away with the physical fluxes in \eqref{eqn.3d-ex}. Likewise, the source terms are truncated away. Therefore, the contributions from the explicit part are zero:
\begin{equation}
\Dr_{i,j,k}^* = (\Delta \vv)_{i,j,k}^* = (\DrE)_{i,j,k}^* = \DB_{i,j,k}^* = 0 \qquad \forall (i,j,k) \in \{1,...,N_x\} \times \{1,...,N_y\} \times \{1,...,N_z\}.
\end{equation}
Consequently, also the deviation in the magnetic energy $\Dm_{i,j,k}^{n+1}$ is zero. 
The constrained transport step after the explicit update does not change anything, because the electric field at the corners in the update \eqref{eqn.Bfc-update} is equal to zero due to the zero velocity in the equilibrium. 
Under these conditions, the elliptic equation for the pressure \eqref{eqn.3d-linSyst_p} simplifies to 
\begin{equation}
\begin{split}
&\frac{p_{i,j,k}^{n+1}}{\gamma-1} \\
&- \frac{\dt^2}{\dx^2} \left( \left( \frac{3}{4} h_{i-1,j,k}^n + \frac{1}{4} h_{i+1,j,k}^n \right) p_{i-1,j,k}^{n+1} - \left( h_{i-1,j,k}^n+h_{i+1,j,k}^n \right) p_{i,j,k}^{n+1} + \left( \frac{1}{4} h_{i-1,j,k}^n+\frac{3}{4} h_{i+1,j,k}^n \right) p_{i+1,j,k}^{n+1} \right) \\
&- \frac{\dt^2}{\dy^2} \left( \left( \frac{3}{4} h_{i,j-1,k}^n + \frac{1}{4} h_{i,j+1,k}^n \right) p_{i,j-1,k}^{n+1} - \left( h_{i,j-1,k}^n+h_{i,j+1,k}^n \right) p_{i,j,k}^{n+1} + \left( \frac{1}{4} h_{i,j-1,k}^n+\frac{3}{4} h_{i,j+1,k}^n \right) p_{i,j+1,k}^{n+1} \right) \\
&- \frac{\dt^2}{\dz^2} \left( \left( \frac{3}{4} h_{i,j,k-1}^n + \frac{1}{4} h_{i,j,k+1}^n \right) p_{i,j,k-1}^{n+1} - \left( h_{i,j,k-1}^n+h_{i,j,k+1}^n \right) p_{i,j,k}^{n+1} + \left( \frac{1}{4} h_{i,j,k-1}^n+\frac{3}{4} h_{i,j,k+1}^n \right) p_{i,j,k+1}^{n+1} \right) \\
=& 
\frac{\pt_{i,j,k}}{\gamma-1} \\
&- \frac{\dt^2}{\dx^2} \left( \left( \frac{3}{4} h_{i-1,j,k}^n + \frac{1}{4} h_{i+1,j,k}^n \right) \pt_{i-1,j,k} - \left( h_{i-1,j,k}^n+h_{i+1,j,k}^n \right) \pt_{i,j,k} + \left( \frac{1}{4} h_{i-1,j,k}^n+\frac{3}{4} h_{i+1,j,k}^n \right) \pt_{i+1,j,k} \right) \\
&- \frac{\dt^2}{\dy^2} \left( \left( \frac{3}{4} h_{i,j-1,k}^n + \frac{1}{4} h_{i,j+1,k}^n \right) \pt_{i,j-1,k} - \left( h_{i,j-1,k}^n+h_{i,j+1,k}^n \right) \pt_{i,j,k} + \left( \frac{1}{4} h_{i,j-1,k}^n+\frac{3}{4} h_{i,j+1,k}^n \right) \pt_{i,j+1,k} \right) \\
&- \frac{\dt^2}{\dz^2} \left( \left( \frac{3}{4} h_{i,j,k-1}^n + \frac{1}{4} h_{i,j,k+1}^n \right) \pt_{i,j,k-1} - \left( h_{i,j,k-1}^n+h_{i,j,k+1}^n \right) \pt_{i,j,k} + \left( \frac{1}{4} h_{i,j,k-1}^n+\frac{3}{4} h_{i,j,k+1}^n \right) \pt_{i,j,k+1} \right)
\end{split}
\end{equation}
which admits the unique solution
\begin{equation}
p_{i,j,k}^{n+1} = \pt_{i,j,k}.
\end{equation}
This means that the right-hand side of the momentum update \eqref{eqn.3d-momentum} becomes zero and consequently the right-hand side of the energy update \eqref{eqn.3d-rhoE} also becomes zero. Thus the updated deviations for all components of the state vector remain zero, which means that the solution $\q$ does not change during one time step, i.e.
\begin{equation}
\q_{i,j,k}^{n+1} = \q_{i,j,k}^n = \qt_{i,j,k} \qquad \forall (i,j,k) \in \{1,...,N_x\} \times \{1,...,N_y\} \times \{1,...,N_z\}.
\end{equation}
\hfill $\square$
\end{proof}

\subsection{Summary of the scheme}
Since various methods are combined in our scheme, we provide an overview of the individual steps of the scheme at this point. 
Let us assume that we know the time-independent equilibrium solution $\qt$ at every point in the computational domain $\Omega$. Then a single time-step in the well-balanced semi-implicit method from time $t=t^n$ to time $t=t^{n+1}$ consists of the following sub-steps:
\begin{enumerate}
\setlength\itemsep{0.5em}
\item Start with the values at the current time level: $\q_{i,j,k}^n$, $\B_{x,i+1/2,j,k}^n$, $\B_{y,i,j+1/2,k}^n$, $\B_{z,i,k,k+1/2}^n$.
\item Compute the deviation: $\Dq_{i,j,k}^n = \q_{i,j,k}^n-\qt_{i,j,k}$.
\item Reconstruct the values at the cell interfaces for the deviation: $\Dq_{i+1/2,j,k}^{L}$, $\Dq_{i+1/2,j,k}^{R}$, $\Dq_{i,j+1/2,k}^{L}$, $\Dq_{i,j+1/2,k}^{R}$, $\Dq_{i,j,k+1/2}^{L}$, $\Dq_{i,j,k+1/2}^{R}$.
\item Perform the explicit update \eqref{eqn.spatial_operators-ex}: $\Dq_{i,j,k}^*$.
\item Start the constrained transport routine: update the staggered magnetic field at the face centers and compute the updated cell centered magnetic field: $\B_{x,i+1/2,j,k}^{n+1}$, $\B_{y,i,j+1/2,k}^{n+1}$, $\B_{z,i,k,k+1/2}^{n+1}$, $\B_{i,j,k}^{n+1}$.
\item Use the computed quantities from the explicit part to solve the elliptic equation for the pressure \eqref{eqn.p} via GMRES: $p_{i,j,k}^{n+1}$.
\item Compute the updated deviation in the momentum: $(\Dru)_{i,j,k}^{n+1}$.
\item Compute the updated deviation in the total energy: $(\DrE)_{i,j,k}^{n+1}$.
\item Recompute the actual solution at the new time level: $\q_{i,j,k}^{n+1}=\Dq_{i,j,k}^{n+1} + \qt_{i,j,k}$.
\setlength\itemsep{0em}
\end{enumerate}
In the case that we do not want to balance a magnetohydrostatic equilibrium, the equilibrium solution $\qt$ can be set to zero.

%--------- SECTION --------------------------------------------------------
\section{Numerical results} \label{sec.numtest} 
%\color{red}scheme settings, cartesian grid, equilibrium solution ... \color{black}

In this section, the second order semi-implicit well-balanced method is applied to various numerical test problems for the Euler and ideal MHD equations. 
%The ratio of specific heat is set to $\gamma=5/3$ if not stated otherwise. 
The time step of the method is calculated in each case based on the CFL condition \eqref{eqn.timestep_im} with $\text{CFL}=0.9$. 
Partially, the results are compared with a non-well-balanced method. This is the same semi-implicit scheme, except for the fact that the equilibrium solution $\qt$ is set to zero. All tests are performed on an uniform Cartesian grid.

%----------
\subsection{Shock tube under gravitational field}
\begin{figure}[!htbp]
\begin{center}
\begin{tabular}{cc}
\includegraphics[width=0.9\textwidth]{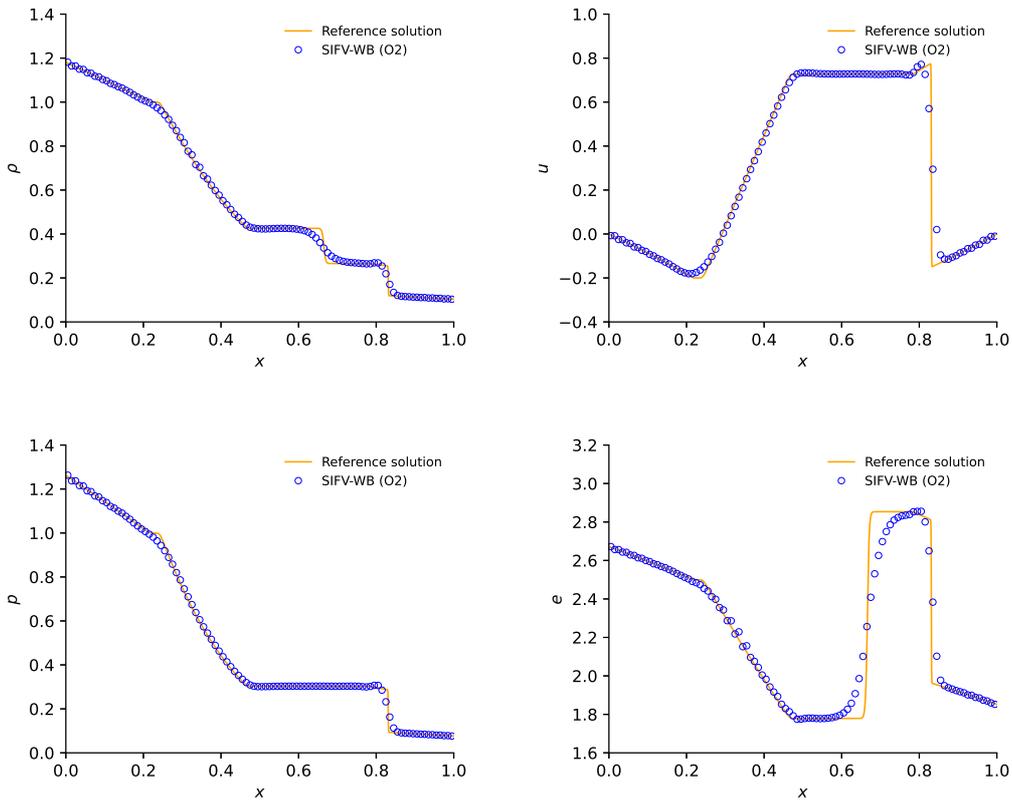} 
\end{tabular}
\caption{Shock tube under gravitational field at time $t_f=0.2$. Comparison against the reference solution for density (top left), horizontal velocity (top right), pressure (bottom left) and internal energy (bottom right).}
\label{fig.Sod_gravity}
\end{center}
\end{figure}

The first test case is the standard Sod shock tube for the one-dimensional Euler equations, to which a gravitational source term with constant acceleration $g_x=-1$ is added \cite{Chandrashekar2015}. The initial conditions in the domain $\Omega=[0,1]$ are given by
\begin{equation}
(\rho,\vv,p) = 
\begin{cases}
(1,\mathbf{0},1), \quad &x< 0.5, \\
(0.125,\mathbf{0},0.1), &x> 0.5,
\end{cases}
\end{equation}
with solid wall boundary conditions. The magnetic field $\B$ is set to zero and the ratio of specific heats is $\gamma=1.4$. The solution at final time $t_f=0.2$ is computed by the second order semi-implicit scheme on 100 cells. Since we have no flow around an equilibrium in this test, we set the equilibrium solution $\qt$ to zero. The numerical solution is compared to a reference solution, which is computed by a fully explicit second order finite volume method on 20000 cells. The results in Fig. \ref{fig.Sod_gravity} show good agreement with the reference solution and are also consistent with solutions in the literature \cite{Chandrashekar2015}. This test demonstrates the capability of the semi-implicit scheme to deal with shocks, thus flows which are not in the low Mach regime.

%\begin{figure}[!htbp]
%\begin{center}
%\begin{tabular}{cc}
%\includegraphics[width=0.47\textwidth]{Figures/Sod_gravity_rho}  &          
%\includegraphics[width=0.47\textwidth]{Figures/Sod_gravity_u}    \\
%\includegraphics[width=0.47\textwidth]{Figures/Sod_gravity_p}  &          
%\includegraphics[width=0.47\textwidth]{Figures/Sod_gravity_e}    \\
%\end{tabular}
%\caption{Shock tube under gravitational field at time $t_f=0.2$. Comparison against the reference solution for density (top left), horizontal velocity (top right), pressure (bottom left) and internal energy (bottom right).}
%\label{fig.Sod_gravity}
%\end{center}
%\end{figure}

\subsection{General Euler Steady State in 1D}
The following test investigates the method's behaviour in the presence of a hydrostatic equilibrium \cite{Desveaux2016}. On the domain $\Omega=[0,1]$ we define the gravitational acceleration in $x$-direction by $g_x = 2 \pi \cos(2\pi x)$ and use the initial values
\begin{equation}
\label{eqn.ic_gss}
(\rho,\vv,p) = \left( 3+2\sin (2\pi x), \ \mathbf{0}, \ 3+3 \sin (2\pi x) - 0.5 \cos (4\pi x) \right).
\end{equation}
The magnetic field $\B$ is set to zero and $\gamma=1.4$. Periodic boundary conditions are imposed on all sides. It can be easily checked that the values in \eqref{eqn.ic_gss} satisfy the equilibrium \eqref{eqn.mhse} and thus should be preserved until any final time $t_f$. Here, we set $t_f=1$. In a first step, we perform this test without defining an equilibrium solution $\qt$. Tab. \ref{tab.GSS} lists the $L_1$ error as well as the experimental order of convergence ($EOC$) for different numbers of grid points $N_x$. Clearly, in this case, the equilibrium is not preserved exactly, but only with a second order of accuracy. In this form, the method is therefore not well-balanced. In a second step we use the initial values as equilibrium solution $\qt$. The results in Tab. \ref{tab.GSS_wb} show that the equilibrium is now maintained exactly. No iterations are needed to solve the linear system for the pressure and the residual is of the size of the machine accuracy.

\begin{table}[!htp]
\begin{center}
\caption{General Steady State for the non well-balanced scheme. $L_1$ error at time $t_f=1$ and the experimental order of convergence for the density, horizontal velocity and pressure using different numbers of cells $N_x$. }
\begin{tabular}{r|cc|cc|cc}
$N_x$ & $L_1(\rho)$ & $EOC(\rho)$ & $L_{1}(u)$ & $EOC(u)$ & $L_1(p)$ & $EOC(p)$ \\
\hline
20    & 2.1451E-003 &       - & 2.4473E-004  &       - & 2.3520E-003 & - \\
40    & 4.2949E-004 & 2.32 & 4.7176E-005  & 2.38 & 4.0100E-004 & 2.55 \\
80    & 1.0332E-004 & 2.19  & 1.3068E-005 & 2.11 & 1.0224E-004 & 2.26 \\
160  & 2.7888E-005 & 2.08 & 3.8141E-006  & 2.00  & 2.8544E-005 & 2.12 \\
\end{tabular}
\label{tab.GSS}
\end{center}
\end{table}

\begin{table}[!htp]
\begin{center}
\caption{General Steady State for the well-balanced scheme. The errors are measured in $L_1$ norm at time $t_f=1$ for the density, horizontal velocity and pressure using different numbers of cells $N_x$. The maximum number of the iterations and the corresponding residuals needed in the linear solver for the pressure system are also reported.}
\begin{tabular}{r|ccc|cc}
$N_x$ & $L_1(\rho)$ & $L_{1}(u)$ & $L_1(p)$ & Iterations & Residual \\
\hline
20    & 0.0000E+00 & 0.0000E+00  & 0.0000E+00 & 0 & 0.0000E+00 \\
40    & 0.0000E+00 & 0.0000E+00  & 0.0000E+00 & 0 & 0.0000E+00 \\
 80    & 0.0000E+00 & 0.0000E+00  & 0.0000E+00 & 0 & 0.0000E+00 \\
160   & 0.0000E+00 & 0.0000E+00  & 0.0000E+00 & 0 & 0.0000E+00 \\
\end{tabular}
\label{tab.GSS_wb}
\end{center}
\end{table}

%----------
\subsection{Isothermal atmosphere for Euler in 2D}
\label{sec.isothermal_equilibrium}
After performing one-dimensional tests, the well-balanced property is now verified for the two-dimensional scheme. For this purpose we set as initial values an isothermal equilibrium \cite{Chandrashekar2015} of the form
\begin{equation}
\label{eqn.ic_isothermal}
(\rho,\vv,p) = \left( 1.21 e^{-1.21(x+y)}, \mathbf{0}, e^{-1.21(x+y)} \right).
\end{equation}
The magnetic field $\B$ is set to zero and the gravitational acceleration is defined by $\g = (-1,-1,0)$. The ratio of specific heats is $\gamma=1.4$. Clearly, this initial data satisfies the equilibrium \eqref{eqn.mhse}. We run this test until a final time $t_f=1$ on a domain $\Omega=[0,1]^2$ with different numbers of cells $(N_x,N_y)$. As boundary conditions we enforce the exact solution. The results in Tab. \ref{tab.IsoAtmosphere} show that indeed no iterations are needed to solve the pressure linear system \eqref{eqn.3d-linSyst_p} because the right hand side is perfectly balanced by the pressure terms. As a consequence, the residual has the order of the machine accuracy and the solution is preserved exactly throughout the simulation.
    
\begin{table}[!htp]
\begin{center}
\caption{Isothermal atmosphere. The errors are measured in $L_1$ norm at time $t_f=1$ for the density, velocity and pressure using different numbers of cells $(N_x,N_y)$. The maximum number of the iterations and the corresponding residuals needed in the linear solver for the pressure system are also reported.}
\begin{tabular}{rr|cccc|cc}
$N_x$ & $N_y$ & $L_1(\rho)$ & $L_{1}(u)$ & $L_{1}(v)$ & $L_1(p)$ & Iterations & Residual \\
\hline
 20 & 20    & 0.0000E+00 & 0.0000E+00 & 0.0000E+00  & 0.0000E+00 & 0 & 0.0000E+00 \\
 40 & 40    & 0.0000E+00 & 0.0000E+00 & 0.0000E+00  & 0.0000E+00 & 0 & 0.0000E+00 \\
 80  & 80   & 0.0000E+00 & 0.0000E+00 & 0.0000E+00  & 0.0000E+00 & 0 & 0.0000E+00 \\
160 & 160 & 0.0000E+00 & 0.0000E+00 & 0.0000E+00  & 0.0000E+00 & 0 & 0.0000E+00 \\
\end{tabular}
\label{tab.IsoAtmosphere}
\end{center}
\end{table}

%----------
\subsection{Perturbation of an isothermal atmosphere for Euler in 2D}
\label{sec.perturbation_isothermal}
\begin{figure}[!htbp]
\begin{center}
\begin{tabular}{cc}
\includegraphics[width=0.9\textwidth]{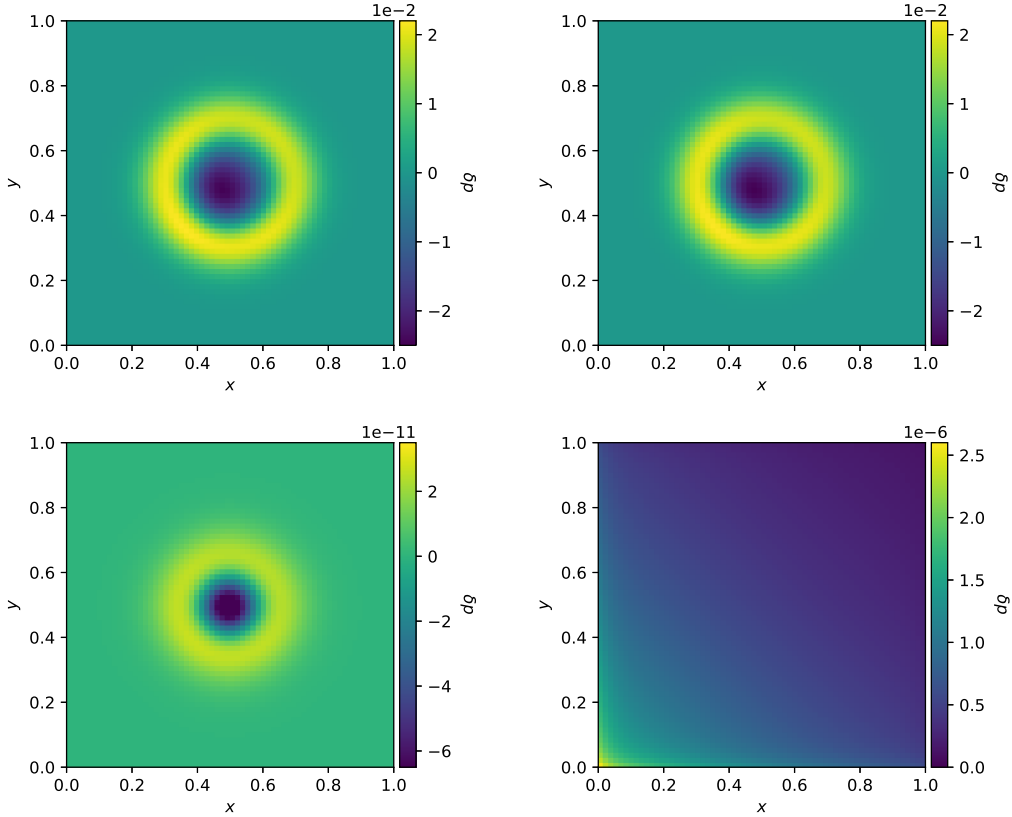} 
\end{tabular}
\caption{Perturbation of an isothermal atmosphere in 2D at time $t=0.15$. Distribution of the pressure perturbation $\delta p$ obtained with an initial strength of $\eta=10^{-1}$ (top row) and $\eta=10^{-10}$ (bottom row) for the well-balanced (left column) and the non-well-balanced (right column) scheme.}
\label{fig.IsoAtmosphere_delta2d}
\end{center}
\end{figure}

The aim of well-balanced methods is not only to maintain hydrostatic equilibria exactly, but in particular to resolve small perturbations of these equilibria. In this sense we now add a perturbation in the pressure to the equilibrium \eqref{eqn.ic_isothermal}, i.e.
\begin{equation}
p = e^{-1.21(x+y)} + \eta e^{-121 \left((x-0.5)^2+(y-0.5)^2 \right)}.
\end{equation}
The parameter $\eta$ regulates the size of the perturbation. We perform the test on a $64\times 64$-grid with the non-well-balanced method as well as with the well-balanced method. For the well-balanced scheme, the unperturbed initial data \eqref{eqn.ic_isothermal} is used for $\qt$. First, we run the test with a rather large perturbation ($\eta = 10^{-1}$). The top row in Fig. \ref{fig.IsoAtmosphere_delta2d} shows the perturbation $\delta p=p(t_f)-p(0)$ at time $t_f=0.15$. For this perturbation, both methods provide a good resolution and perform qualitatively equally well. In a second simulation, we choose a smaller perturbation ($\eta = 10^{-10}$). In this case, the non-well-balanced method is no longer able to resolve the perturbation, while the well-balanced method resolves the small perturbation as good as the large one before. These results are consistent with the results from the literature \cite{Chandrashekar2015,Semplice2019,Birke2023}.

%\begin{figure}[!htbp]
%	\begin{center}
%		\begin{tabular}{cc}
%			\includegraphics[width=0.47\textwidth]{Figures/2d_Euler_equilibrium_delta_1e-1_WB}  &          
%			\includegraphics[width=0.47\textwidth]{Figures/2d_Euler_equilibrium_delta_1e-1_noWB}    \\
%			\includegraphics[width=0.47\textwidth]{Figures/2d_Euler_equilibrium_delta_1e-10_WB}  &          
%			\includegraphics[width=0.47\textwidth]{Figures/2d_Euler_equilibrium_delta_1e-10_noWB}    \\
%		\end{tabular}
%		\caption{Perturbation of an isothermal atmosphere in 2D at time $t=0.15$. Distribution of the pressure perturbation $\delta p$ obtained with an initial strength of $\eta=10^{-1}$ (top row) and $\eta=10^{-10}$ (bottom row) for the well-balanced (left column) and the non-well-balanced (right column) scheme.}
%		\label{fig.IsoAtmosphere_delta2d}
%	\end{center}
%\end{figure}

\subsection{MHD Steady State in 2D}
\label{sec.mhd_equilibrium}
The following test is intended to show that the method can also maintain equilibria exactly, which have a nonzero magnetic field. For this purpose, we modify the equilibrium from Sect. \ref{sec.isothermal_equilibrium} by adding a magnetic field. The initial values on the domain $\Omega = [0,1]^2$ are then given by
\begin{equation}
\label{eqn.ic_mhd_equilibrium}
(\rho, \vv, p, B_x, B_y, B_z ) = 
\left( 2.21 e^{-(x+y)}, \ \mathbf{0}, \ 1.21 e^{-(x+y)}, \ e^{-\frac{1}{2}(x+y)}, \ -e^{-\frac{1}{2}(x+y)}, \ 0 \right),
\end{equation}
and we set $\gamma=1.4$. The boundary conditions rely on the exact solution. The $L_1$-errors for the well-balanced scheme at time $t_f=1$ are given in Tab. \ref{tab.mhd_equilibrium} and show that the equilibrium is preserved exactly. Consequently, the method is also well-balanced for MHD equilibria as stated in Theorem \ref{th.wb_pressure}.

\begin{table}[!htp]
\begin{center}
\caption{MHD equilibrium. The errors are measured in $L_1$ norm at time $t_f=1$ for the density, velocity and pressure using different numbers of cells $N_x$. The maximum number of the iterations and the corresponding residuals needed in the linear solver for the pressure system are also reported.}
\begin{tabular}{rr|cccc|cc}
$N_x$ & $N_y$ & $L_1(\rho)$ & $L_{1}(u)$ & $L_{1}(v)$ & $L_1(p)$ & Iterations & Residual \\
\hline
 20 & 20    & 0.0000E+00 & 0.0000E+00 & 0.0000E+00  & 0.0000E+00 & 0 & 0.0000E+00 \\
 40 & 40    & 0.0000E+00 & 0.0000E+00 & 0.0000E+00  & 0.0000E+00 & 0 & 0.0000E+00 \\
 80  & 80   & 0.0000E+00 & 0.0000E+00 & 0.0000E+00  & 0.0000E+00 & 0 & 0.0000E+00 \\
160 & 160 & 0.0000E+00 & 0.0000E+00 & 0.0000E+00  & 0.0000E+00 & 0 & 0.0000E+00 \\
\end{tabular}
\label{tab.mhd_equilibrium}
\end{center}
\end{table}

\subsection{Perturbation of the MHD Steady State in 2D}
\label{sec.mhd_perturbation}
\begin{figure}[!htbp]
\begin{center}
\begin{tabular}{cc}
\includegraphics[width=0.9\textwidth]{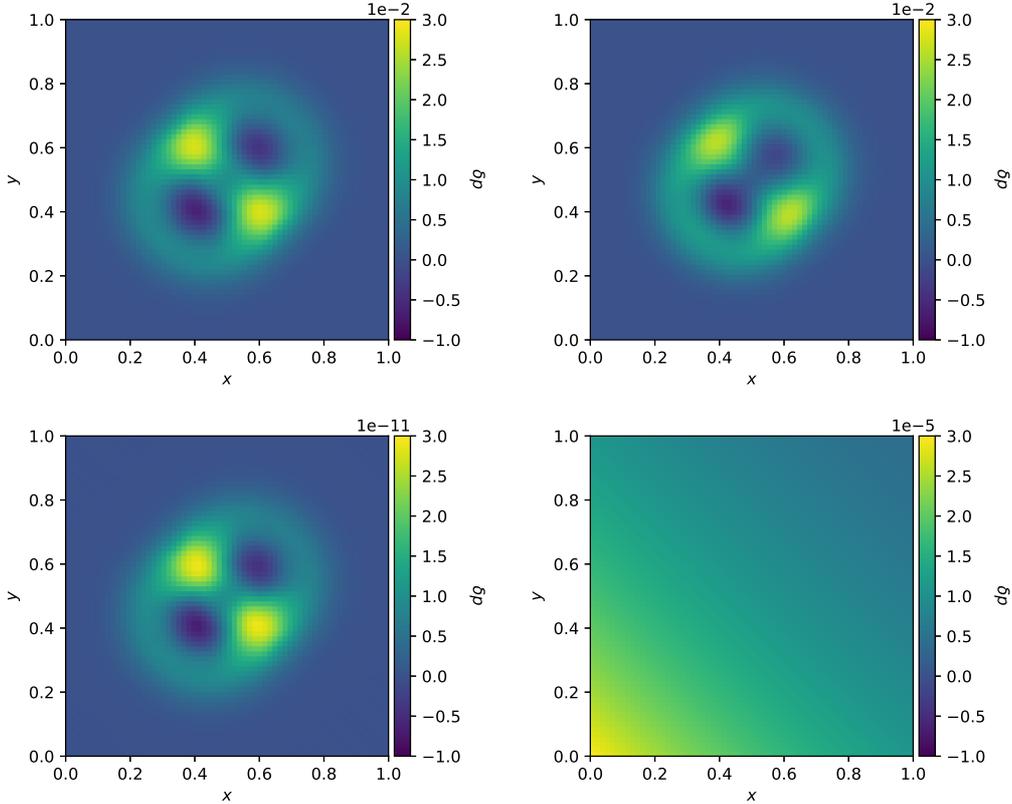} 
\end{tabular}
\caption{Perturbation of a MHD equilibrium in 2D at time $t=0.15$. Distribution of the pressure perturbation $\delta p$ obtained with an initial strength of $\eta=10^{-1}$ (top row) and $\eta=10^{-10}$ (bottom row) for the well-balanced (left column) and the non-well-balanced (right column) scheme.}
\label{fig.mhd_equilibrium_delta2d}
\end{center}
\end{figure}

As already done for the isothermal Euler equilibrium, we also add a perturbation to the MHD equilibrium  \eqref{eqn.ic_mhd_equilibrium}. The initial pressure is then given by
\begin{equation}
p = 1.21 e^{-(x+y)} + \eta e^{-100 \left((x-0.5)^2+(y-0.5)^2 \right)}.
\end{equation}
This test is also carried out with a large ($\eta=10^{-1}$) and a small ($\eta=10^{-10}$) perturbation. Fig. \ref{fig.mhd_equilibrium_delta2d} shows the results for the well-balanced and the non-well-balanced method at final time $t_f=0.15$. The behavior is similar to the one for the Euler equilibrium in Sect. \ref{sec.perturbation_isothermal}. The non-well-balanced method can only resolve the large perturbation, whereas the well-balanced method accurately reproduces both perturbations. It should be noted that the altered shape of the Gaussian perturbation in comparison to the one of the Euler equilibrium results from the newly introduced force of the magnetic field.

\subsection{Euler vortex in a gravitational field}
\label{sec.gravity_vortex}
\begin{figure}[!htbp]
\begin{center}
\begin{tabular}{cc}
\includegraphics[width=0.9\textwidth]{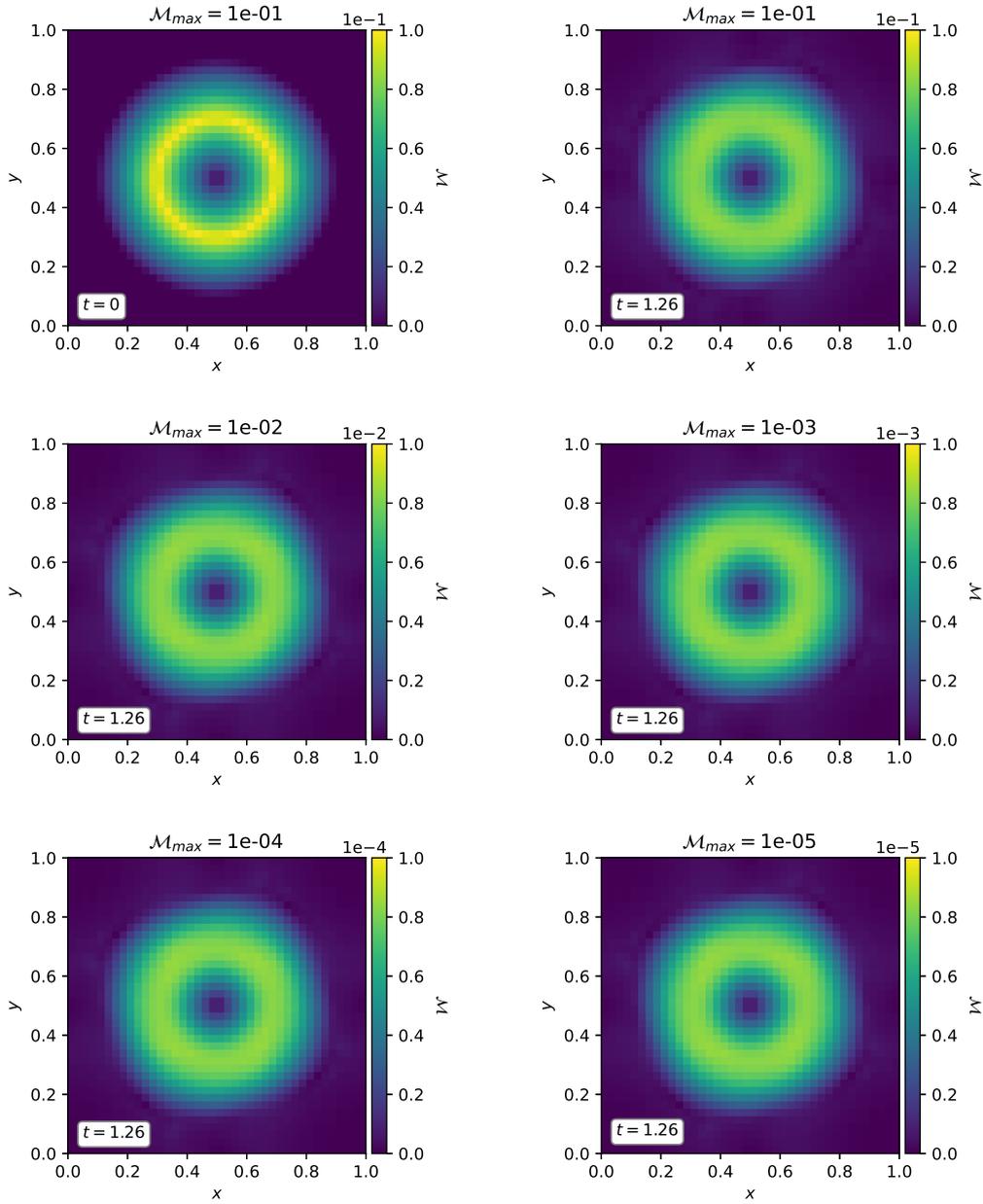} 
\end{tabular}
\caption{Distribution of the local Mach number $\mathcal{M}$ at $t=0$ (top left) and after one turn for different maximum Mach numbers $\mathcal{M}_{max}$.}
\label{fig.gravity_vortex}
\end{center}
\end{figure}

\begin{figure}[!htbp]
\begin{center}
\begin{tabular}{cc}
\includegraphics[width=0.7\textwidth]{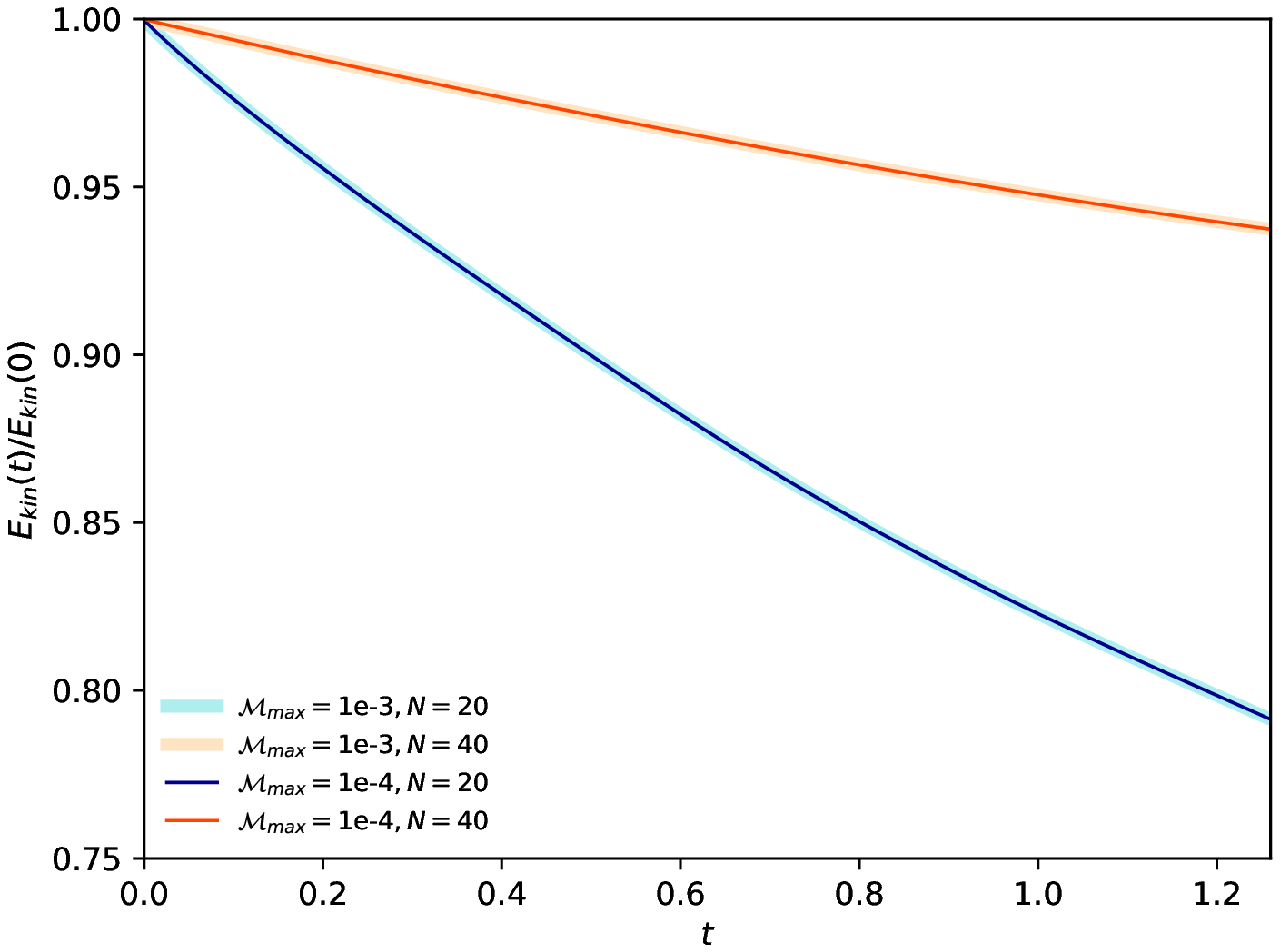} 
\end{tabular}
\caption{Time evolution of the total kinetic energy for different maximum Mach numbers $\mathcal{M}_{max}$ and grid resolutions $N=N_x=N_y$.}
\label{fig.ekin}
\end{center}
\end{figure}

After demonstrating the well-balanced property in the former tests, this section is devoted to show the low Mach property of the scheme. For this purpose, we use a vortex for the Euler equations that is characterised by the balance of the pressure gradient with the centrifugal and the gravitational forces, which makes the vortex axisymmetric, stationary and having zero radial velocity. The vortex is placed on top of a hydrostatic equilibrium solution, 
%\begin{equation}
%\rho_{hs} = \exp \left( - \frac{\Phi}{\mathcal{M}_{max}^2} \right), \quad \vv_{hs} = 0, \quad p_{hs}= \mathcal{M}_{max}^2 \rho_{hs},
%\end{equation}
which is why $\qt$ is set to this equilibrium. The detailed initial data can be found in \cite{Thomann2020}. \\
%\rev{Should we give the complete initial data? Quite long ...} \\
We compute the solution after one full turn of the vortex ($t_f=1.26$) for different initial maximum Mach numbers $\mathcal{M}_{max}$ on a $40 \times 40$-grid. Fig. \ref{fig.gravity_vortex} shows the distributions of the local Mach number $\mathcal{M}$ at time $t_f$ and additionally the initial distribution for $\mathcal{M}_{max}=10^{-1}$ (top left). It becomes clear that the numerical scheme adds a certain amount of numerical dissipation to the approximate solution in comparison to the initial data, which dampens the local Mach number. The vortex, however, is well resolved in each case and there is no qualitative difference between the solutions for different $\mathcal{M}_{max}$. This indicates that the numerical dissipation is independent of the Mach number thanks to the implicit discretisation of the acoustic terms in the flux splitting.\\
This finding is underlined by the time evolution of the total kinetic energy. Since the domain is a closed setup due to the periodic boundary conditions, the total kinetic energy should be conserved over time in the low Mach limit. Therefore, the loss of kinetic energy is a good measure for the amount of numerical dissipation. Fig. \ref{fig.ekin} shows the time evolution of the kinetic energy for different $\mathcal{M}_{max}$ and grid resolutions $N=N_x=N_y$. It turns out that the loss of kinetic energy only depends on the grid resolution and not on the Mach number.
The results are consistent with those in \cite{Thomann2020}.

%\subsection{Hot Bubble}

%--------- SECTION --------------------------------------------------------
\section{Conclusions} 
\label{sec.concl}

In this work we have presented a well-balanced semi-implicit scheme for the ideal MHD equations with gravitational source terms. The method is based on splitting the equations into the convective part, which is discretized explicitly in time, and the acoustic part, which is discretized implicitly. An elliptic equation for the pressure is obtained, which is solved in the implicit part without adding any numerical dissipation. This makes the method particularly suitable for problems in the low Mach regime. Numerical results show indeed that the numerical dissipation is independent of the Mach number. At the same time, the explicit finite volume method for the nonlinear convective terms ensures stability in the presence of shocks. A constrained transport method is embedded in the numerical scheme to ensure the solenoidal property of the magnetic field at the discrete level.\\
In order to maintain magnetohydrostatic equilibria exactly and to handle small perturbations of these equilibria even on coarse grids, the method utilizes the deviation well-balancing technique. The resulting scheme solves the equations for the deviation of the solution from an \textit{a priori} known equilibrium solution. We have proven that, in the equilibrium, the pressure equation has only one solution, which is the equilibrium pressure itself. As a result, the solution is preserved exactly in the case of a magnetohydrostatic equilibrium. Thus, the scheme is well-balanced. Numerical tests demonstrate that arbitrary equilibria are preserved exactly and even very small perturbations are resolved accurately. \\
Future developments will concern the extension of the method to unstructured meshes along the lines of \cite{sifvdg,fvvem}, in order to address more challenging and physically relevant problems in plasma flows as tokamak simulations \cite{WBtokamak_Fambri}.

%--------- SECTION --------------------------------------------------------
\section*{Acknowledgments}

CB acknowledges the support by the German Research Foundation (DFG) under the project No. KL 566/22-1.	
WB received financial support by Fondazione Cariplo and Fondazione CDP (Italy) under the project No. 2022-1895 and by the Italian Ministry of University and Research (MUR) with the PRIN Project 2022 No. 2022N9BM3N. WB is member of the GNCS-INdAM (\textit{Istituto Nazionale di Alta Matematica}) group.

	\section*{Declarations}
	
	\paragraph{Funding.} CB received financial support by the German Research Foundation (DFG) under the project No. KL 566/22-1. WB received financial support by Fondazione Cariplo and Fondazione CDP (Italy) under the project No. 2022-1895 and by the Italian Ministry of University (MUR) with the PRIN Project 2022 No. 2022N9BM3N.
	
	\paragraph{Conflicts of interest.} The authors declare that they have no conflict of interest.
	
	\paragraph{Availability of data and material.} Data and material are available upon reasonable request addressed to the corresponding author.
	
	\paragraph{Code availability.} The code is written in Fortran programming language and is available upon reasonable request addressed to the corresponding author.
	
	\paragraph{Ethics approval.} Not applicable.
	
	\paragraph{Consent to participate.} Not applicable.
		
	\paragraph{Consent for publication.} Not applicable.
	
	\section*{Data availability}
	The datasets generated during the current study are available from the corresponding author upon reasonable request.

	\bibliographystyle{plain}
	\bibliography{biblio}   % name your BibTeX data base
	
\end{document}